\documentclass[10pt]{article}

\usepackage{amssymb}

\usepackage[utf8]{inputenc}
\usepackage{graphicx}

\usepackage{nicematrix}
\usepackage{colortbl}
\definecolor{Gray}{gray}{0.9}
\usepackage{ragged2e}
\newcolumntype{g}{>{\columncolor{Gray}}r}
\newcolumntype{h}{>{\columncolor{Gray}}r}

\usepackage{arxiv}
\usepackage{relsize}
\usepackage[capitalize]{cleveref}
\usepackage{algorithm}
\usepackage{algpseudocode}
\usepackage{ifthen}
\usepackage[font=normalsize]{subcaption}
\usepackage{color}
\usepackage{wrapfig}
\usepackage{float}
\usepackage{bm}
\usepackage{multirow}
\usepackage{array} 
\graphicspath{{images/}}

\usepackage{mathtools}

\newcommand{\cred}[1]{{\color{black}  #1}}
\newcommand{\credd}[1]{{\color{black}  #1}}
\newcommand{\crep}[1]{{\color{black}  #1}}
\newcommand{\cblu}[1]{{\color{black}  #1}}

\newcommand{\p}{_{N}}
\newcommand{\q}{\zeta}

\newtheorem{theorem}{Theorem}[section]

\newtheorem{proposition}[theorem]{Proposition}
\newtheorem{definition}[theorem]{Definition}
\newtheorem{remark}[theorem]{Remark}

\usepackage{amsopn}
\DeclareMathOperator{\diag}{diag}

\makeatletter
\newcommand{\srcsize}{\@setfontsize{\srcsize}{5pt}{5pt}}
\makeatother

\newcommand{\abs}[1]{\left\lvert#1\right\rvert}
\newcommand{\norm}[1]{\left\lVert#1\right\rVert}

\renewcommand{\O}{\mathrm{O}}

\newcommand{\OO}[1]{\ifthenelse{\equal{#1}{0}}{\mathrm{O}(h)}{
		\ifthenelse{\equal{#1}{+1}}{\mathrm{O}(h)}{
			\ifthenelse{\equal{#1}{-1}}{\mathrm{O}h)}{
				\mathrm{O}(h)
}}}}  

\newcommand{\g}{g'(\hat{x}_i)}
\newcommand{\gi}{g'_i}

\newcommand{\numberset}{\mathbb}

\newcommand{\R}{\numberset{R}}

\newcommand{\ord}{t}

\title{Multigrid for two-sided fractional differential equations discretized by finite volume elements on graded meshes}

\author{
Marco Donatelli\\
Department of Science and High Technology\\
Insubria University\\
marco.donatelli@uninsubria.it\\
\And
Rolf Krause\\
Faculty of Informatics\\
University of Italian Switzerland\\
rolf.krause@usi.ch\\
\And
Mariarosa Mazza\\
Department of Science and High Technology\\
Insubria University\\
mariarosa.mazza@uninsubria.it\\
\And
Ken Trotti\\
Faculty of Informatics\\
University of Italian Switzerland\\
ken.trotti@usi.ch
}

\begin{document}
	\maketitle
	
\begin{abstract}
It is known that the solution of a conservative steady-state \cred{two-sided} fractional diffusion problem can exhibit singularities near the boundaries. As consequence of this, and due to the conservative nature of the problem, we adopt a finite volume elements discretization approach over a generic non-uniform mesh. We focus on grids mapped by a smooth function which consist in a combination of a graded mesh near the singularity and a uniform mesh where the solution is smooth. Such a choice gives rise to Toeplitz-like discretization matrices and thus allows a low computational cost of the matrix-vector product and a detailed spectral analysis. The obtained spectral information is used to develop an ad-hoc parameter free multigrid preconditioner for GMRES, which is numerically shown to yield good convergence results in presence of graded meshes mapped by power functions that accumulate points near the singularity. The approximation order of the considered graded meshes is numerically compared with the one of a certain composite mesh given in literature that still leads to Toeplitz-like linear systems and is then still well-suited for our multigrid method. Several numerical tests confirm that power graded meshes result in lower approximation errors than composite ones and that our solver has a wide range of applicability.
\end{abstract}

\keywords{
two-sided fractional problems \and finite volume elements methods \and Toeplitz matrices \and spectral distribution \and multigrid methods
}



\section{Introduction}\label{sec_intro}
We consider a conservative steady-state \cred{two-sided} Fractional Diffusion Equation (FDE) of order $2-\beta$, $0<\beta<1$, with inhomogeneous Dirichlet boundary-value conditions 
\cite{Conservative_caputo,Conservative_caputo2}, i.e.,
\begin{equation}\label{eq_FDE}
\left\{
\begin{split}
&-\frac{\mathrm{d}}{\mathrm{d}x}\left(K(x)\big(\gamma\ _0\mathcal{D}_x^{1-\beta}+(1-\gamma)\  _x\mathcal{D}_1^{1-\beta}\big)\right)u(x)=f(x),\qquad 0<x<1,\\
&\ u(0)=u_l,\quad u(1)=u_r,
\end{split}
\right.
\end{equation}
where $K(x)$ is a positive diffusion coefficient, $f(x)$ is the source term, $u_l,u_r$ are the Dirichlet boundary values and $0\leq\gamma\leq 1$ indicates the anisotropy in the diffusion, i.e., $\gamma\approx 0$ and $\gamma\approx 1$ imply a strong forward and backward diffusivity, respectively. By   $_0\mathcal{D}_x^{1-\beta},_x\mathcal{D}_1^{1-\beta}$ we denote the left and right Caputo fractional derivatives (to be defined in the next section), while the fractional derivative operators $\frac{\mathrm{d}}{\mathrm{d}x}\ _0\mathcal{D}_x^{1-\beta}$ and $\frac{\mathrm{d}}{\mathrm{d}x}\ _x\mathcal{D}_1^{1-\beta}$ are known as Riemann–Liouville–Caputo fractional derivatives \cite{RLC} or as Patie–Simon fractional derivatives \cite{PS1,PS2,PS3}.

This kind of FDEs may exhibit singularities near the boundaries even in case of smooth coefficients; see, e.g., \cite{Kopteva,stynes} where the Caputo derivative is in time or \cite{wang} where the Caputo derivative is in space. As a consequence, non-uniform mesh-based discretization methods should be used. In this light, and due to the conservative nature of the problem, in \cite{wang} the authors propose a Finite Volume Element (FVE) discretization approach over a composite mesh made of a graded part and a uniform one.

In this paper, we adopt the FVE approach over a generic non-uniform mesh and we explicitly provide all the coefficients of the resulting linear system. In particular, we focus on graded meshes mapped by power functions near the singularity. As for the case of the composite mesh used in \cite{wang}, also for meshes mapped by functions, the FVE discretization leads to structured linear systems with a Toeplitz-like pattern. Moreover, the regularity of the mapping function allows a spectral study of the coefficient matrices. Extending the work done in \cite{FVE_RL_marco} in presence of uniform meshes, here we perform a spectral analysis of the coefficient matrices combining the generating function of the Toeplitz part with the first derivative of the function that defines the graded mesh. We use the retrieved information to design an ad-hoc multigrid preconditioner, which is parameter free since the relaxing parameter of the Jacobi smoother is estimated through the approach introduced in \cite{paper_MGS}. 

A wide number of numerical tests compares our proposal applied to both power graded meshes and the composite mesh given in \cite{wang} with the circulant preconditioner proposed therein. Concerning the parameter that defines the power mapping, we resort to the results in \cite{Kopteva}, where the authors deal with a time-fractional diffusion equation, with Caputo derivative in time of order $\alpha\in(0,1)$, discretized through an $L1$ approximation scheme. \cred{We numerically show that such a choice fits also within our setting and that the resulting convergence order is similar to the one obtained in \cite{Kopteva}, while lower approximation errors are obtained with respect to the composite meshes used in \cite{wang}. Finally, compared to the circulant preconditioner proposed in \cite{wang} for composite meshes, our solver demonstrates linearly convergent over both graded and composite meshes.}

The paper is organized as follows. In Section \ref{sec:preliminaries}, we define the fractional derivative operators and recall a few results on Toeplitz and generalized locally Toeplitz sequences. In Section \ref{sec_FVscheme}, we provide the full FVE discretization over arbitrary meshes, then we focus on meshes mapped by power functions and in Section \ref{sec:spectral} we retrieve the spectral information of the involved discretization matrices. Our multigrid proposal is discussed in Section \ref{sec:vcycle} and numerically tested in Section \ref{sec_results}. Finally, in Section \ref{sec:conclusion} we draw conclusions.

\section{Preliminaries}\label{sec:preliminaries}
This section contains various preliminaries on fractional derivatives (Section~\ref{sub:frac}) and Toeplitz matrices (Section~\ref{sec:toeplitz}) needed in the rest of the paper. In Section~\ref{sec:toeplitz}, we also briefly introduce the Generalized Locally Toeplitz (GLT) theory which extends the spectral results for symmetric Toeplitz matrices to more general cases.
\subsection{Fractional derivatives}\label{sub:frac}
For a given function with absolutely integrable first derivative on $[0,1]$, the right-handed and left-handed Caputo fractional derivatives of order $1-\beta$, $0<\beta<1$ are defined by 
\begin{align*}
_0\mathcal{D}_x^{1-\beta}g(x)&\coloneqq \frac{1}{\Gamma(\beta)}\int_0^x (x-s)^{\beta-1}\frac{\mathrm{d}g(s)}{\mathrm{d}s}\mathrm{d}s,\qquad
_x\mathcal{D}_1^{1-\beta}g(x)\coloneqq -\frac{1}{\Gamma(\beta)}\int_x^1 (s-x)^{\beta-1}\frac{\mathrm{d}g(s)}{\mathrm{d}s}\mathrm{d}s,
\end{align*}
with $\Gamma(\cdot)$ the Euler Gamma function. Another common definition of fractional derivatives that asks for less regularity on the function is due to Riemann and Liouville
\begin{align*}
\ _0\mathcal{R}_x^{1-\beta}g(x)&\coloneqq \frac{1}{\Gamma(\beta)}\frac{\mathrm{d}}{\mathrm{d}s}\int_0^x (x-s)^{\beta-1}g(s)\mathrm{d}s,\qquad
\ _x\mathcal{R}_1^{1-\beta}g(x)\coloneqq -\frac{1}{\Gamma(\beta)}\frac{\mathrm{d}}{\mathrm{d}s}\int_x^1 (s-x)^{\beta-1}g(s)\mathrm{d}s.
\end{align*}
In this case it is enough that $g$ is absolutely continuous. The Riemann-Liouville derivatives relate to the Caputo ones as follows (see Equations (2.4.8)-(2.4.9) at page 91 of \cite{north_holland})
\begin{align}\label{eq:rel}
\begin{split}
\ _0\mathcal{R}_x^{1-\beta}g(x)&=\ _0\mathcal{D}_x^{1-\beta}g(x)
+\frac{x^{\beta-1}}{\Gamma(\beta)}g(0),\qquad
\ _x\mathcal{R}_1^{1-\beta}g(x)=\ _x\mathcal{D}_1^{1-\beta}g(x)
+\frac{(1-x)^{\beta-1}}{\Gamma(\beta)}g(1).
\end{split}
\end{align}
Therefore, in case of homogeneous Dirichlet boundary conditions the two definitions coincide. 
%
\subsection{Toeplitz and generalized locally Toeplitz sequences}\label{sec:toeplitz}

In this subsection, we recall the definition of Toeplitz sequences generated by a function and we recall some key properties of the GLT class which will be used in Section~\ref{sec:spectral} to provide spectral information on the coefficient matrices resulting after a certain non-uniform FVE discretization of equation \eqref{eq_FDE}.
\begin{definition}\label{def_toeplitz}
	A \emph{Toeplitz matrix} $T_N\in\mathbb{C}^{N\times N}$ has constant coefficients along the diagonals, namely
	$\left[T_N\right]_{i,j}=t_{i-j},\ i,j=1,...,N$.
	If $\lbrace t_k\rbrace_{k\in\mathbb{Z}}$ are the Fourier coefficients of a function $f\in L^1([-\pi,\pi])$, i.e.,
	$t_k=\frac{1}{2\pi}{\int_{0}^{2\pi}}f(\theta)\text{e}^{-\mathrm{i}k\theta}\mathrm{d}\theta,$
	the function $f$ is called the \emph{generating function} of $\{T_N\}_N$, and we write $T_N=T_N(f)$.
\end{definition}
The GLT class is a matrix-sequence algebra obtained as a closure under some algebraic operations between Toeplitz and diagonal matrix-sequences generated by functions (to be defined below). It includes matrix-sequences coming from the discretization of differential operators with various techniques, such as finite differences, finite elements, Isogeometric Analysis, etc. The formal definition of GLT class is difficult and involves a heavy notation, therefore in the following we just introduce those among its properties that we need for our studies (for a more detailed discussion see \cite{serra_vol1}).

\begin{definition}\label{def_symbol_diag}
	A matrix-sequence whose $N$-th element is a diagonal matrix $D_N=[d_{i,j}]_{i,j=1}^N\in\mathbb{R}^{N\times N}$ such that $\ d_{i,i}=d\left(\frac{i}{N}\right),\ i=1,...,N$, with $d\!:[0,1]\rightarrow\mathbb{C}\ $ a Riemann-integrable function, is called \emph{diagonal sampling sequence}.
\end{definition}
The functions $f$ in Definition \ref{def_toeplitz} and $d$ in Definition \ref{def_symbol_diag} allow to estimate the spectrum of the matrix-sequences $\{T_N(f)\}_N$ and $\{D_N\}$, respectively, in the following sense.
\begin{definition}\label{def_distribution_eig}
	Let $f:G\rightarrow\mathbb{C}$ be a measurable function, defined on a measurable set $G\subset\mathbb{R}^k$ with $k\geq 1$, $0<m_k(G)<\infty$, where $m_k(G)$ is the Lebesgue measure of the set $G$. Let $C_0(\mathbb{K})$ be the set of continuous functions with compact support over $\mathbb{K}\in\lbrace\mathbb{R}_0^+,\mathbb{C}\rbrace$ and let $\lbrace A_N\rbrace_{N}$ be a sequence of matrices of size $N$ with eigenvalues $\lambda_j(A_N), j = 1,...,N$ and singular values $\sigma_j(A_N), j = 1,...,N$.
	\begin{itemize}
		\item $\lbrace A_N\rbrace_{N}$ is distributed as the pair $(f,G)$ in the sense of the eigenvalues, in {formulae} $\lbrace A_N\rbrace_{N}\!\sim_\lambda\!(f,G)$, if the following limit relation holds for all $F\!\in\! C_0(\mathbb{C})$ 
		\begin{equation}\label{distribution:eig}
		\lim\limits_{N\rightarrow\infty}\frac{1}{N}\sum\limits_{j=1}^{N}F(\lambda_j(A_N))=\frac{1}{m_k(G)}\int_G F(f(t))\mathrm{d}t.
		\end{equation}
		{In this case, we refer to the function $f$ as (spectral) symbol.}
		\item $\lbrace A_N\rbrace_{N}$ is distributed as the pair $(f,G)$ in the sense of the singular values, in {formulae} $\lbrace A_N\rbrace_{N}\!\sim_\sigma\!(f,G)$, if the following limit relation holds for all $F\!\in\! C_0(\mathbb{R}_0^+)$
		\begin{equation}\label{distribution:sv}
		\lim\limits_{N\rightarrow\infty}\frac{1}{N}\sum\limits_{j=1}^{N}F(\sigma_j(A_N))=\frac{1}{m_k(G)}\int_G F(\abs{f(t)})\mathrm{d}t.
		\end{equation}
		{In this case, we refer to the function $f$ as singular value symbol.}
	\end{itemize}
\end{definition} 
\begin{remark} \label{rem:informal}
	An informal interpretation of the limit relation \eqref{distribution:eig} {(resp. \eqref{distribution:sv})} is that when $N$ is sufficiently large, the eigenvalues {(resp. singular values)} of $A_N$ can be approximated by a sampling of $f$ {(resp. $\abs{f}$)} on a uniform mesh over the set $G$, up to a relatively small number of potential outliers and where ``relatively small" means $o(N)$.
\end{remark}

Throughout, we use the following notation $$\lbrace A_N\rbrace_{N}\sim_{GLT} {\psi}(x,\theta),\qquad (x,\theta)\in [0,1]\times[-\pi,\pi],$$ to say that the sequence $\lbrace A_N\rbrace_{N}$ is a GLT sequence with symbol $\psi(x,\theta)$.

Here we report {five} main features of the GLT class.
\begin{description}
	\item[GLT1] Let $\{A_N\}_{N}\sim_{\rm GLT}\psi(x,\theta)$ with $\psi:G\rightarrow \mathbb{C}$, $G=[0,1]\times[-\pi,\pi]$, then $\{A_N\}_{N}\sim_\sigma(\psi,G)$. If the matrices $A_N$ are Hermitian, then  $\{A_N\}_{N}\sim_\lambda(\psi,G)$.
	\item[GLT2] The set of GLT sequences form a $*$-algebra, i.e., it is closed under linear combinations, products, and transposed conjugation. Moreover, it is closed under inversion whenever the symbol vanishes, at most, in a set of zero Lebesgue measure. Hence, the sequence obtained via algebraic operations on a finite set of input GLT sequences is still a GLT sequence and its symbol is obtained by following the same algebraic manipulations on the corresponding symbols of the input GLT sequences.
	\item[GLT3] Every Toeplitz sequence $\{T_N(f)\}\p$ generated by a $L^1([-\pi,\pi])$ function $f(\theta)$ is such that $\{T_N(f)\}\p\sim_{\rm GLT}f(\theta)$, with the specifications reported in item {\bf GLT1}. Every diagonal sampling sequence $\{D_{N}(a)\}\p$, where $a(x)$ is a Riemann integrable function {in [0,1]}, is such that $\{D_N(a)\}\p\sim_{\rm GLT}a(x)$.
	\item[GLT4] Every sequence distributed as the constant zero in the singular value sense is a GLT sequence with symbol zero, and viceversa. In formulae, $\{A_N\}\p\sim_\sigma(0,G)$, $G=[0,1]\times[-\pi,\pi]$, if and only if $\{A_N\}\p\sim_{\rm GLT}0$.
	\item[GLT5] Let $\lbrace A_N\rbrace_{N}\sim_{\mathrm{GLT}} \psi(x,\theta)$, $G=[0,1]\times[-\pi,\pi]$. If we assume that 
	$$\mathrm{lim}_{N\to\infty}\frac{\norm{A_N-A_N^{\mathrm{H}}}_{\mathrm{tr}}}{N}{=0},$$
	where $\norm{\cdot}_\mathrm{tr}$ is the trace norm, i.e., the sum of the singular values, then $\psi(x,\theta)$ is necessarily a real-valued function and $\lbrace A_N\rbrace_{N}\sim_{\lambda} \psi(x,\theta)$.
\end{description}
%
The first GLT result that we need in the next sections is reported in Proposition \ref{prop_distribution} and concerns the symbol of a diagonal-times-Toeplitz matrix-sequence.
\begin{proposition}[\cite{serra_vol1}]\label{prop_distribution}
	Let $\lbrace D_N\rbrace_{N}$ be a sequence of diagonal sampling matrices with symbol $d:[0,1]\rightarrow\mathbb{R}_{>0}\ $, and $\lbrace T_N{(f)}\rbrace_{N}$ be a sequence of Hermitian Toeplitz matrices with symbol $f:[-\pi,\pi]\rightarrow\R$, then
	$$\lbrace D_N T_N(f)\rbrace_{N}\sim_\lambda \big(d(x)f(\theta),[0,1]\times[-\pi,\pi]\big).$$
\end{proposition}
Notice that Proposition \ref{prop_distribution} is a consequence of the similitude transformation $D_N^{-\frac{1}{2}} D_N T_N(f) D_N^{\frac{1}{2}} = D_N^{\frac{1}{2}}T_N(f) D_N^{\frac{1}{2}}$ and of the hermitianity of $D_N^{\frac{1}{2}}T_N(f) D_N^{\frac{1}{2}}$ which, from \textbf{GLT1-3}, ensure
\begin{equation*}
\lbrace D_N^{\frac{1}{2}}T_N(f) D_N^{\frac{1}{2}}\rbrace_{N} \sim_\lambda \big(\sqrt{d(x)}f(\theta)\sqrt{d(x)}=d(x)f(\theta),\ [0,1]\times[-\pi,\pi]\big).
\end{equation*}
In case the diagonal-times-Toeplitz structure is hidden, one can resort to the notion of approximating class of sequences and to the GLT result reported in Theorem \ref{th_simbolo_acs}, which allows to find the symbol of a \lq difficult' matrix-sequence by means of \lq simpler' matrix-sequences.
\begin{definition}\label{def_acs}
	Let $\lbrace A_N\rbrace_{N}$ be a matrix-sequence and let $\lbrace\lbrace B_{N,M}\rbrace_{N}\rbrace_M$ be a sequence of matrix-sequences. We say that $\lbrace\lbrace B_{N,M}\rbrace_{N}\rbrace_M$ is an \emph{approximating class of sequences} (a.c.s.) for $\lbrace A_N\rbrace_{N}$ if the following condition is met: $\forall M\ \exists\ N_M$ such that for $N\geq N_M$,
	\begin{equation*}
	A_N=B_{N,M}+R_{N,M}+N_{N,M},\quad \mathrm{rank}(R_{N,M})\leq c(M)N,\quad \norm{N_{N,M}}_2\leq \omega(M),
	\end{equation*}
	where $\norm{\cdot}_2$ is the spectral norm and $N_M,c(M),\omega(M)$ depend only on $M$ with
	\begin{equation*}
	\lim\limits_{M\rightarrow\infty} c(M)=\lim\limits_{M\rightarrow\infty}\omega(M)=0.
	\end{equation*}
\end{definition}
\begin{theorem}[\cite{serra_vol1}]\label{th_simbolo_acs}
	Let $\lbrace A_N\rbrace_{N}$ be a matrix-sequence. If there exists an a.c.s. $\lbrace\lbrace B_{N,M}\rbrace_{N}\rbrace_M$ for $\lbrace A_N\rbrace_{N}$ such that $\lbrace\lbrace B_{N,M}\rbrace_{N}\rbrace_M\sim_\sigma (f_{N,M},G)$, with $f_{N,M}$ that converges in measure to $f$, then
	$$\lbrace A_N\rbrace_{N}\sim_\sigma (f,G).$$
\end{theorem}

\section{Finite volume elements scheme}\label{sec_FVscheme}
The FVE approach consists in restricting the admissible solutions $u(x)$ {of the FDE in \eqref{eq_FDE}} to a certain finite element space, partitioning the definition interval $[0,1]$ as $\cup_{i=1}^n \Omega_i$, where $\mu(\Omega_i\cap \Omega_j)\!=\!0,\ i\!\neq\!j$, with $\mu$ the Lebesgue measure, and finally integrating equation \eqref{eq_FDE} over $\Omega_i$. In our specific case we consider the unknown $u(x)$ to belong to the space of the piecewise linear polynomial functions. In the following, we provide a full discretization over a generic non-structured mesh (Section \ref{sec:generic-non-uniform-mesh}), then we consider the special case of a uniform mesh (Section \ref{sec:uniform-mesh}) and make a comparison with the discretization in \cite{FVE_RL_marco}. Finally, in Section \ref{sec_grid}, we {focus on} two non-uniform structured meshes, {obtained as a combination of a non-uniform part and a uniform one}.


\subsection{Generic non-uniform mesh}\label{sec:generic-non-uniform-mesh}
Let $N\in\mathbb{N}$ and denote by $\lbrace x_i\rbrace_{i=0}^{N+1}$ a generic mesh on $[0,1]$, such that $x_i>x_{i-1},\ \forall i=1,...,N+1$ with $x_0=0, x_{N+1}=1$, then {define}
\begin{equation*}
\tilde{u}(x)=\sum_{i=1}^{N}u_{i}\phi_i(x)+u_l\phi_0(x)+u_r\phi_{N+1}(x),
\end{equation*}
where $\lbrace \phi_i\rbrace_{i=0}^{N+1}$ is the set of hat (linear) functions with
\begin{align*}
\phi_i(x)&=
\begin{cases}
\frac{x-x_{i-1}}{h_i}, &x\in(x_{i-1},x_{i})\\
\frac{x_{i+1}-x}{h_{i+1}}, &x\in(x_i,x_{i+1})\\
0, & \text{otherwise}
\end{cases}\qquad \text{for}\ i=1,...,N,
\end{align*}
{and
\begin{align*}
\phi_0(x)&=
\begin{cases}
	\frac{x_{1}-x}{h_{1}}, &x\in(x_0,x_{1})\\
	0, & \text{otherwise}
\end{cases},\qquad 
\phi_{N+1}(x)=
\begin{cases}
	\frac{x-x_{N}}{h_{N+1}}, &x\in(x_{N},x_{N+1})\\
	0, & \text{otherwise}
\end{cases},
\end{align*}
where} $h_i=x_i-x_{i-1},\ i=1,...,N+1$ {is} the step length. {Replacing} $u(x)$ with $\tilde{u}(x)$ {in equation~\eqref{eq_FDE} and} integrating over $\Omega_i=[x_{i-\frac{1}{2}},x_{i+\frac{1}{2}}]$, where $x_{i-\frac{1}{2}}=\frac{x_i+x_{i-1}}{2}$, equation~\eqref{eq_FDE} can be written as the linear system
\begin{equation}\label{eq_LINEAR_system_discretizedeq}
A_N u=b,
\end{equation}
where $b\in\R^{N}$ and $A_N \in\R^{N\times N}$, with
\begin{align*}
a_{i,j}&=-K(x)\left(\gamma\ _0\mathcal{D}_x^{1-\beta}+(1-\gamma)\ _x\mathcal{D}_1^{1-\beta}\right)\phi_j(x)\Big\vert_{x=x_{i-\frac{1}{2}}}^{x=x_{i+\frac{1}{2}}},\\
b_i&=\int_{x_{i-\frac{1}{2}}}^{x_{i+\frac{1}{2}}}f(x)\mathrm{d}x+K(x)\left(\gamma\ _0\mathcal{D}_x^{1-\beta}+(1-\gamma)\ _x\mathcal{D}_1^{1-\beta}\right)(u_l\phi_0(x)+u_r\phi_{N+1}(x))\Big\vert_{x=x_{i-\frac{1}{2}}}^{x=x_{i+\frac{1}{2}}},
\end{align*}
{for $i,j=1,...,N$. Explicitly{, let $K_i=K(x_i)$ and $f_i=f(x_i)$, then} the entries of $b$ are}
{\setlength{\abovedisplayskip}{2pt}\begin{flalign*}
	\quad b_1=&f_1\frac{h_1+h_{2}}{2}+
	\frac{K_{\frac{3}{2}}}{\Gamma(\beta+1)}\Big[\frac{u_L\gamma}{h_1}\left((x_{\frac{3}{2}}-x_1)^\beta-x_{\frac{3}{2}}^\beta\right) +\frac{u_R(1-\gamma)}{h_{N+1}}\left((1-x_{\frac{3}{2}})^\beta-(x_N-x_{\frac{3}{2}})^\beta\right)\Big]+&\\
	&-\frac{K_{\frac{1}{2}}}{\Gamma(\beta+1)}\Big[
	-\frac{u_L\gamma}{h_1}x_{\frac{1}{2}}^\beta+
	(1-\gamma)\left(-\frac{u_L}{h_1}(x_1-x_\frac{1}{2})^{\beta}+\frac{u_R}{h_{N+1}}\left((1-x_{\frac{1}{2}})^\beta-(x_N-x_{\frac{1}{2}})^\beta\right)\right)
	\Big],
	\end{flalign*}}
{\setlength{\abovedisplayskip}{2pt}\begin{flalign*}
	\quad b_i=&f_i\frac{h_i+h_{i+1}}{2}+
	\frac{K_{i+\frac{1}{2}}}{\Gamma(\beta+1)}\Big[
	\frac{u_L\gamma}{h_1}\left((x_{i+\frac{1}{2}}-x_1)^\beta-x_{i+\frac{1}{2}}^\beta\right)+\frac{u_R(1-\gamma)}{h_{N+1}}\left((1-x_{i+\frac{1}{2}})^\beta-(x_N-x_{i+\frac{1}{2}})^\beta\right)\Big]+&\\
	&-\frac{K_{i-\frac{1}{2}}}{\Gamma(\beta+1)}\Big[
	\frac{u_L\gamma}{h_1}\left((x_{i-\frac{1}{2}}-x_1)^\beta-x_{i-\frac{1}{2}}^\beta\right)+\frac{u_R(1-\gamma)}{h_{N+1}}\left((1-x_{i-\frac{1}{2}})^\beta-(x_N-x_{i-\frac{1}{2}})^\beta\right)\Big],
	\end{flalign*}}
{for $i=2,...,N-1$, and}
{\setlength{\abovedisplayskip}{2pt}\begin{flalign*}
	\quad b_N=&f_N\frac{h_N+h_{N+1}}{2}+
	\frac{K_{N+\frac{1}{2}}}{\Gamma(\beta+1)}\Big[
	\frac{u_L\gamma}{h_1}\left((x_{N+\frac{1}{2}}-x_1)^\beta-x_{N+\frac{1}{2}}^\beta\right)+ \frac{ u_R\gamma}{h_{N+1}}(x_{N+\frac{1}{2}}-x_N)^\beta+\frac{u_R(1-\gamma)}{h_{N+1}}(1-x_{N+\frac{1}{2}})^\beta\Big]+&\\
	&-\frac{K_{N-\frac{1}{2}}}{\Gamma(\beta+1)}\Big[
	\frac{u_L\gamma}{h_1}\left((x_{N-\frac{1}{2}}-x_1)^\beta-x_{N-\frac{1}{2}}^\beta\right)+\frac{u_R(1-\gamma)}{h_{N+1}}\left((1-x_{N-\frac{1}{2}})^\beta-(x_N-x_{N-\frac{1}{2}})^\beta\right)\Big].
	\end{flalign*}}
{The entries of $A_N$, for $i=1,...,N$, are}
{\setlength{\abovedisplayskip}{2pt}\begin{flalign*}
	\qquad a_{i,i-k}=&\frac{K_{i-\frac{1}{2}}}{\Gamma(\beta+1)}\gamma\left[
	\frac{(\frac{h_i}{2}+\sum_{j=1}^{k} h_{i-j})^\beta-(\frac{h_i}{2}+\sum_{j=1}^{k-1} h_{i-j})^\beta}{h_{i-k}}+
	\frac{(\frac{h_i}{2}+\sum_{j=1}^{k-2} h_{i-j})^\beta-(\frac{h_i}{2}+\sum_{j=1}^{k-1} h_{i-j})^\beta}{h_{i-k+1}}
	\right]+&\\
	&-\frac{K_{i+\frac{1}{2}}}{\Gamma(\beta+1)}\gamma\left[
	\frac{(\frac{h_{i+1}}{2}+\sum_{j=0}^k h_{i-j})^\beta-(\frac{h_{i+1}}{2}+\sum_{j=0}^{k-1} h_{i-j})^\beta}{h_{i-k}}+
	\frac{(\frac{h_{i+1}}{2}+\sum_{j=0}^{k-2} h_{i-j})^\beta-(\frac{h_{i+1}}{2}+\sum_{j=0}^{k-1} h_{i-j})^\beta}{h_{i-k+1}}
	\right],
	\end{flalign*}}
for $2\leq k\leq i-1$, and
{\setlength{\abovedisplayskip}{2pt}\begin{flalign*}
	\qquad a_{i,i-1}=&\frac{K_{i-\frac{1}{2}}}{\Gamma(\beta+1)}\left[
	\gamma\frac{(h_{i-1}+\frac{h_i}{2})^\beta-(\frac{h_{i}}{2})^\beta}{h_{i-1}}-\frac{(\frac{h_i}{2})^\beta}{h_i}
	\right]+\nonumber\\
	&-\frac{K_{i+\frac{1}{2}}}{\Gamma(\beta+1)}\gamma\left[
	\frac{(h_{i-1}+h_i+\frac{h_{i+1}}{2})^\beta-(h_i+\frac{h_{i+1}}{2})^\beta}{h_{i-1}}+\frac{(\frac{h_{i+1}}{2})^\beta-(h_i+\frac{h_{i+1}}{2})^\beta}{h_{i}}
	\right]
	;&\\
	a_{i,i}=&\frac{K_{i-\frac{1}{2}}}{\Gamma(\beta+1)}\left[
	\frac{(\frac{h_i}{2})^\beta}{h_i}+(1-\gamma)\frac{(\frac{h_{i}}{2})^\beta-(h_{i+1}+\frac{h_{i}}{2})^\beta}{h_{i+1}}
	\right]-\frac{K_{i+\frac{1}{2}}}{\Gamma(\beta+1)}\left[
	\gamma\frac{(h_i+\frac{h_{i+1}}{2})^\beta-(\frac{h_{i+1}}{2})^\beta}{h_i}-\frac{(\frac{h_{i+1}}{2})^\beta}{h_{i+1}}
	\right]
	;&\\
	a_{i,i+1}=&\frac{K_{i-\frac{1}{2}}}{\Gamma(\beta+1)}(1-\gamma)\left[
	\frac{(h_{i+1}+\frac{h_i}{2})^\beta-(\frac{h_i}{2})^\beta}{h_{i+1}}+\frac{(h_{i+1}+\frac{h_{i}}{2})^\beta-(h_{i+2}+h_{i+1}+\frac{h_{i}}{2})^\beta}{h_{i+2}}
	\right]+&\\
	&-\frac{K_{i+\frac{1}{2}}}{\Gamma(\beta+1)}\left[
	\frac{(\frac{h_{i+1}}{2})^\beta}{h_{i+1}}+(1-\gamma)\frac{(\frac{h_{i+1}}{2})^\beta-(h_{i+2}+\frac{h_{i+1}}{2})^\beta}{h_{i+2}}
	\right],
	\end{flalign*}}
and finally,
{\setlength{\abovedisplayskip}{2pt}\begin{flalign*}
	a_{i,i+k}=&\frac{K_{i-\frac{1}{2}}}{\Gamma(\beta+1)}(1-\gamma)\left[
	\frac{(\frac{h_i}{2}+\sum_{j=1}^k h_{i+j})^\beta-(\frac{h_i}{2}+\sum_{j=1}^{k-1} h_{i+j})^\beta}{h_{i+k}}+
	\frac{(\frac{h_i}{2}+\sum_{j=1}^{k} h_{i+j})^\beta-(\frac{h_i}{2}+\sum_{j=1}^{k+1} h_{i+j})^\beta}{h_{i+k+1}}
	\right]+&\\
	&-\frac{K_{i+\frac{1}{2}}}{\Gamma(\beta+1)}(1-\gamma)\left[
	\frac{(\frac{h_{i+1}}{2}+\sum_{j=2}^k h_{i+j})^\beta-(\frac{h_{i+1}}{2}+\sum_{j=2}^{k-1} h_{i+j})^\beta}{h_{i+k}}+
	\frac{(\frac{h_{i+1}}{2}+\sum_{j=2}^{k} h_{i+j})^\beta-(\frac{h_{i+1}}{2}+\sum_{j=2}^{k+1} h_{i+j})^\beta}{h_{i+k+1}}
	\right],
	\end{flalign*}}
for $2\leq k\leq N-i$.

\begin{remark}\label{rem_b0_b1}
	\cred{Let us consider $K(x)=1$.} When $\beta=0$, we have $a_{i,i-k}=a_{i,i+k}=0,\forall k\geq 2$ and the dense structure of $A_N$ collapses into the tridiagonal matrix representing the 1D discrete \cred{Laplacian operator}, which does not depend on $\gamma$ anymore. \cred{On the contrary, when $\beta=1$ we still have $a_{i,i-k}=a_{i,i+k}=0,\forall k\geq 2$ independently of $\gamma$ and $A_N$ becomes a skew-symmetric matrix.}
\end{remark}

\subsection{Uniform mesh}\label{sec:uniform-mesh}
Under the conditions
\begin{equation}\label{eq_Conditions}
K(x)=K,\qquad \gamma=\frac{1}{2},\qquad h_i=h=\frac{1}{N+1},\ i=1,...,N+1,
\end{equation}
with $c=\frac{Kh^{\beta-1}}{2^{\beta}\Gamma(\beta+1)}$, it holds
\begin{equation}             
\begin{split}\label{eq_matrix_uniform_grid}
a_{i,i-k}=&\frac{c}{2}\left[
3(2k+1)^\beta-3(2k-1)^\beta+(2k-3)^\beta-(2k+3)^\beta
\right],\qquad 2\leq k\leq i-1;\\
a_{i,i-1}=&\frac{c}{2}\left[
3^{\beta+1}-4-5^\beta
\right];\\
a_{i,i}=&\frac{c}{2}\left[
6-2\cdot 3^\beta
\right];\\
a_{i,i+1}=&\frac{c}{2}\left[
3^{\beta+1}-4-5^\beta
\right];\\
a_{i,i+k}=&\frac{c}{2}\left[
3(2k+1)^\beta-3(2k-1)^\beta+(2k-3)^\beta-(2k+3)^\beta
\right],\qquad 2\leq k\leq N-i.
\end{split}
\end{equation}
Therefore, under the assumptions in equation \eqref{eq_Conditions}, matrix $A_N$ in \eqref{eq_LINEAR_system_discretizedeq} is a symmetric Toeplitz matrix and coincides with the coefficient matrix considered in \cite{Wang-Du}, where the authors consider a FVE discretization of \eqref{eq_FDE} with Riemann-Liouville fractional derivative operators in place of Caputo's. This is indeed not surprising since, from \eqref{eq:rel} and by $\phi_i(0)= \phi_i(1)=0,\ \forall i=1,...,N$,  we have 
\begin{align*}
\ _0\mathcal{R}_x^{1-\beta}\phi_i(x)=\ _0\mathcal{D}_x^{1-\beta}\phi_i(x),\\
\ _x\mathcal{R}_1^{1-\beta}\phi_i(x)=\ _x\mathcal{D}_1^{1-\beta}\phi_i(x),
\end{align*}
which means that the only difference between the FVE discretization of equation \eqref{eq_FDE} on uniform meshes and the discretized equation in \cite{Wang-Du} lies in the right-hand side.

\subsection{Graded and composite meshes}\label{sec_grid}

The discretization of equation \eqref{eq_FDE} over uniform meshes yields matrices with a Toeplitz structure, which allows fast matrix-vector product in $\mathrm{O}(N\log N)$, while in case of a generic non-uniform mesh discretization the Toeplitz structure is lost. On the other hand, the solution of the FDE in \eqref{eq_FDE} may exhibit singularities near the boundaries, therefore uniform grids should be avoided and non-uniform meshes should be preferred.

In order to deal with the singularity and at the same time \credd{to do} not completely lose the structure of the coefficient matrices, in the following we consider {two mixed approaches} of graded mesh near the singularity and a uniform mesh where the solution is smooth. This yields matrices with a partial Toeplitz structure that can be exploited to allow a fast matrix-vector product. {For the sake of simplicity,} in the following we only consider singularities at $x=0$, however, the approach can be straightforwardly extended to the case of singularities at $x=1$ or at both boundaries. 

\paragraph{Graded meshes} 
{Let $N\in\mathbb{N}$ and consider the uniform grid $\lbrace \hat{x}_i\rbrace_{i=0}^{N+1}$ with $\hat{x}_i=ih$, $i=0,...,N+1$, and $h=\tfrac{1}{N+1}$,
then
\begin{equation}\label{eq_graded_mesh}
	\lbrace {x}_i\rbrace_{i=0}^{N+1},\qquad x_i=g(\hat{x}_i),
\end{equation}
is the non-uniform grid generated by projection of the uniform mesh $\lbrace \hat{x}_i\rbrace_{i=0}^{N+1}$ through the endomorphism $g(x):[0,1]\to[0,1]$.
We consider endomorphisms of the form}
\begin{equation}\label{eq_grid_g}
g_{q,\mathbf{\epsilon}}(x)=
\begin{cases}
x^q, & 0\leq x \leq \epsilon_1,\\
ax^2+bx+c,& \epsilon_1 \leq x \leq \epsilon_1+\epsilon_2,\\
mx+p,& \epsilon_1+\epsilon_2\leq x\leq 1,
\end{cases}
\end{equation}
with ${\bf\epsilon}=(\epsilon_1,\epsilon_2)$, $0<\epsilon_1+\epsilon_2\leq1,\ {\epsilon_2>0}$, and $a,b,c,m,p$ such that $g_{q,\mathbf{\epsilon}}\in \mathrm{C}^1([0,1])$. Proposition \ref{prop_existence_grid}, proved in appendix \ref{appendix_2}, shows that $g_{q,\mathbf{\epsilon}}(x)$ is well-defined{, i.e., given $q,\epsilon_1,\epsilon_2$ \credd{there exist unique} $a,b,c,m,p$ such that $g_{q,\mathbf{\epsilon}}\in \mathrm{C}^1([0,1])$}.
{
	\begin{proposition}\label{prop_existence_grid}
		Let $g_{q,\mathbf{\epsilon}}(x)$ be as in \eqref{eq_grid_g}, with ${\bf\epsilon}=(\epsilon_1,\epsilon_2)$. Then, for $0<\epsilon_1+\epsilon_2\leq1$ with ${\epsilon_2>0}$, function $g_{q,\mathbf{\epsilon}}(x)$ is well-defined.
	\end{proposition}
}
{Note that in the case where $\epsilon_2=0$, the interval $[\epsilon_1,\epsilon_1+\epsilon_2]=\{\epsilon_1\}$ and the quadratic function in \eqref{eq_grid_g} disappears, leading to a loss in smoothness of $g_{q,\epsilon}$.} \credd{Indeed, the function $g_{q,\epsilon}$ has been defined such that in $[0,\epsilon_1]$ accumulates grid points near the singularity at the origin, in $[\epsilon_1+\epsilon_2,1]$ gives a uniform mesh, and in $[\epsilon_1,\epsilon_1+\epsilon_2]$ is a quadratic function that acts as a smooth connection of length $\epsilon_2$ between the singular part and the uniform mesh, \credd{with the only purpose to} increase the smoothness of the whole function $g_{q,\mathbf{\epsilon}}$. Therefore, it is clear that $\epsilon_1+\epsilon_2$ represents the length of the non-uniform part of the grid over the interval $[0,1]$.}

{When $\mathbf{\epsilon}$ is fixed, the only free parameter of $g_{q,\epsilon}$ is $q$ which we choose to be}
{
	\begin{equation}\label{eq_q}
	q=q_\beta=\frac{1+\beta}{1-\beta},
	\end{equation}
}
as done in \cite{Kopteva}, where a Caputo time-fractional derivative is involved and a $L1$ approximation is considered. Therein the authors proved the convergence order to be $1+\beta$ over quasi-graded meshes, which are meshes asymptotically close to graded \credd{ones} mapped by $g_{q,\mathbf{\epsilon}}(x)$ with $\mathbf{\epsilon}=(1,0)$.\\
When $N$ is large, $g_{q,\mathbf{\epsilon}}(x)$ could map a grid that has too short intervals, i.e, $h_1=g_{q,\epsilon}(\hat{x}_1)\ll 10^{-16}$. \credd{Therefore 
we replace $q_\beta$ with $q=\frac{\log(10^{-16})}{-\log(N+1)}$ such that}
\begin{equation}\label{eq_lower_cap}
	h_1 = 10^{-16}.
\end{equation}

\paragraph{Composite mesh} For our numerical comparisons, we will consider also the composite mesh used in \cite{wang}, which has been proved to be effective in the case where $\beta\approx 1$. Let $N_1,N_2\in\mathbb{N}$ and consider an uniform mesh with step $h=\frac{1}{\cblu{N_2+1}}$.\\
Then we divide the interval $[0,h]$ into $N_1+1$ subintervals, whose lengths from left to right are $h_i,\ i=1,...,N_1\cblu{+1}$ with
	\begin{equation*}
	\begin{cases}
	h_i=2^{-N_1}h,     &\text{if}\ i=1;\\
	h_i=2^{i-2-N_1}h, &\text{if}\ i=2,...,N_1\cblu{+1}.
	\end{cases}
	\end{equation*}
	Therefore, the grid points are 
	{
	\begin{equation}\label{eq_composite_mesh_x}
	\begin{cases}
		x_i=0,& \text{if}\ i=0\\
		x_i=2^{i-1-N_1}h,     &\text{if}\ i=1,...,N_1;\\
		x_{i+N_1}=ih, &\text{if}\ i=1,...,N_2;\\
		x_i=1,& \text{if}\ i=N_1+N_2+1.
	\end{cases}
	\end{equation}
	We fix $N\in\mathbb{N}$ and then we choose
	\begin{equation}\label{eq_composite_mesh}
	N_1=g_G(N)\quad \text{and}\quad N_2=N-N_1,
	\end{equation}
	such that the total amount of grid points is $N$, where 
	$g_G(N):\mathbb{N}\rightarrow\mathbb{N}$, e.g., $g_G(N)=\lfloor\sqrt{N}\rfloor$ or $g_G(N)=\lfloor \log_2{N}\rfloor$, with $\lfloor \cdot\rfloor$ being the floor function. \\
	Note that in \cite{wang} the authors consider $h=\frac{1}{N_2}$, leading to coefficient matrices of size $N-1$, and they also first fix $N_2$ and then choose $N_1\approx \sqrt{N_2}$. Our choices are needed to have matrices of size $N$ as with graded meshes and hence to provide a meaningful comparison between the two approaches.}
	



\section{Spectral properties of the coefficient matrices}\label{sec:spectral}
In this section, we first recall the spectral symbol of the coefficient matrices $A_N$ in presence of uniform meshes already given in \cite{FVE_RL_marco}. Then, following the idea in \cite{serra_vol1} (p. 212, section 10.5.4), we compute the symbol of the coefficient matrices when considering non-uniform grids mapped by functions. In both cases we fix $\gamma=\frac{1}{2}$.

\begin{figure}
	\centering
	\begin{subfigure}[b]{0.48\textwidth}
		\centering
		\includegraphics[width=\textwidth]{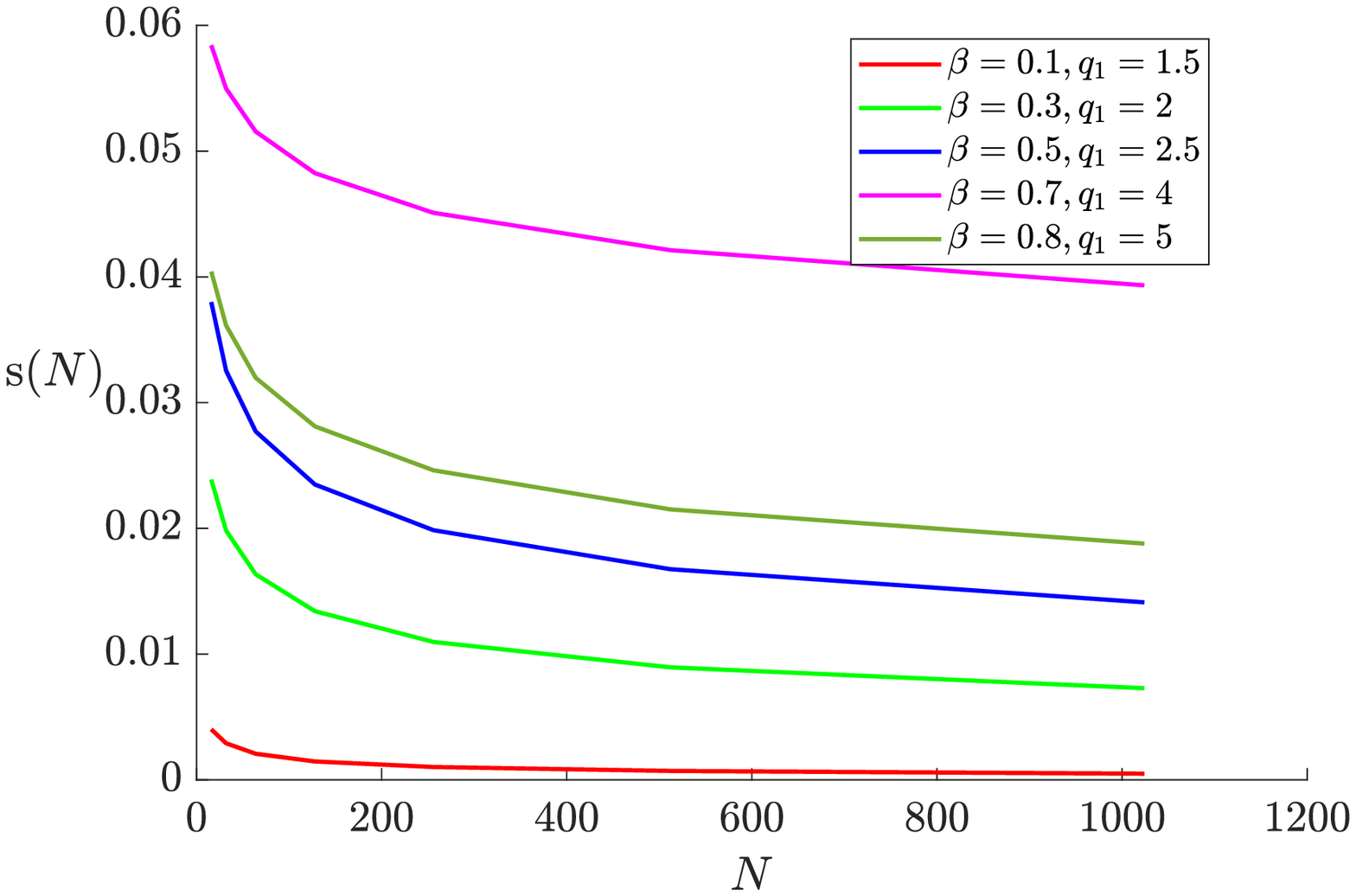}
		\caption{Values of $s(N)$ varying $N$, $\beta\in[0.1,0.8]$ and $q_1<\frac{2-\beta}{1-\beta}$.}
		\label{fig_q1}
	\end{subfigure}
	\hfill
	\begin{subfigure}[b]{0.48\textwidth}
		\centering
		\includegraphics[width=\textwidth]{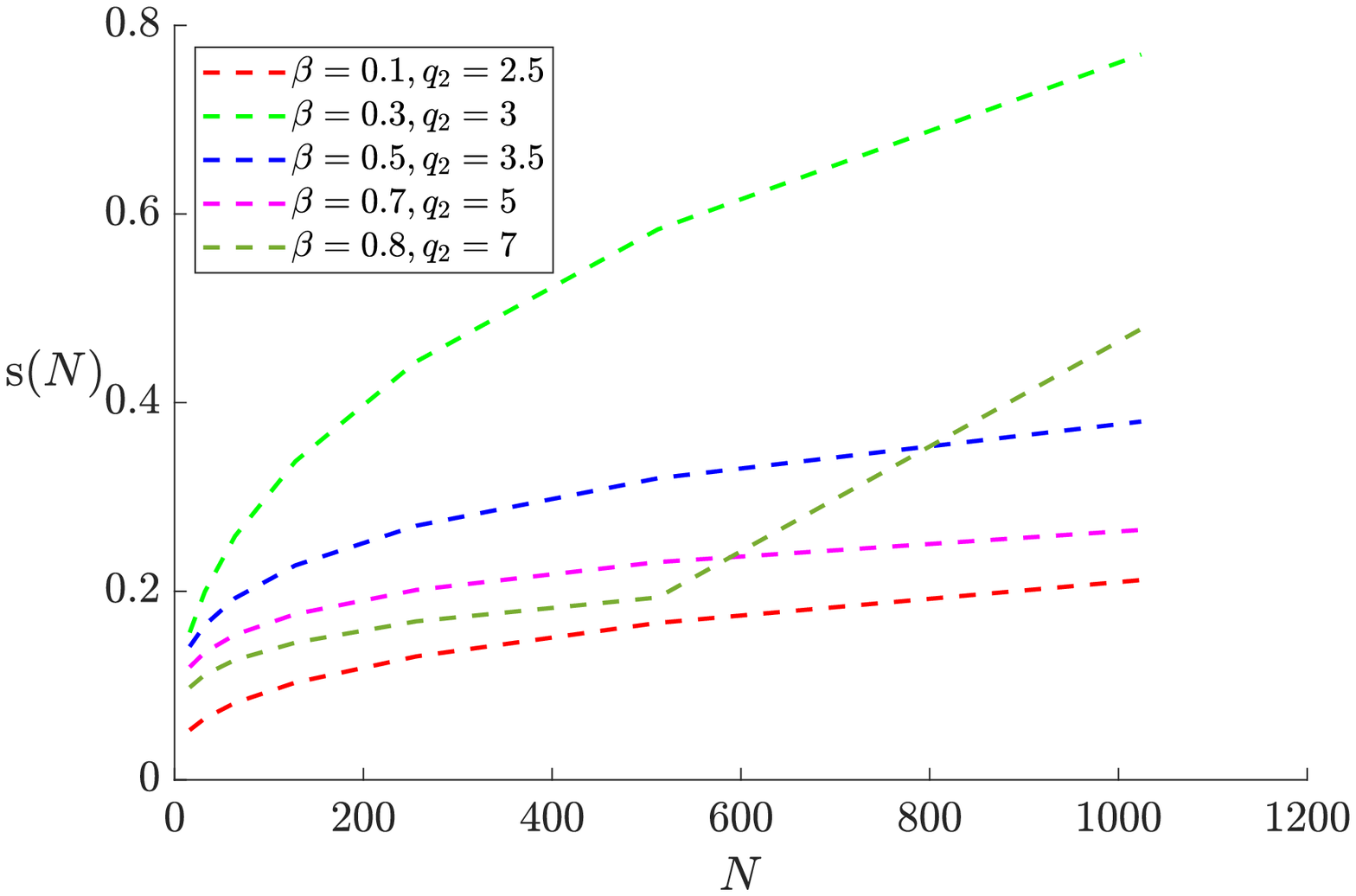}
		\caption{Values of $s(N)$ varying $N$, $\beta\in[0.1,0.8]$ and $q_2>\frac{2-\beta}{1-\beta}$.}
		\label{fig_q2}
	\end{subfigure}
	\caption{Plot of the sequence $s(N)$ varying $N$ for different combinations of $\beta$ and $q$.}\label{fig_sequence}
\end{figure}

\begin{figure}
	\centering
	\includegraphics[width=0.6\textwidth]{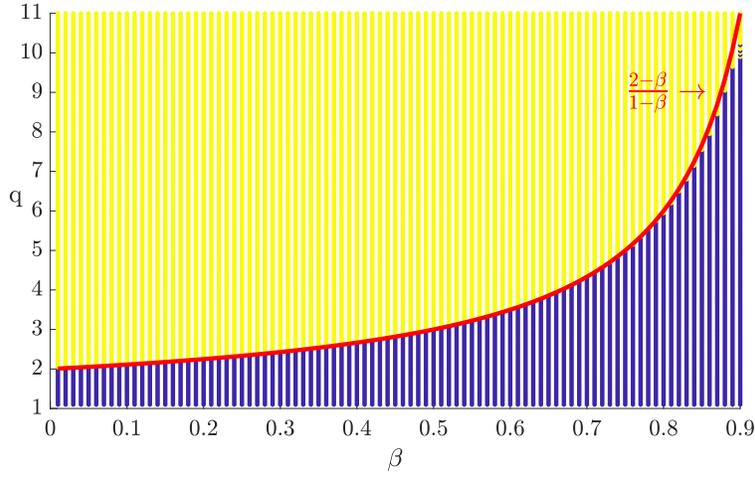}
	\caption{Region where \textbf{GLT5} holds true (in blue), varying $\beta,q$.}
	\label{fig_GLT5_q}
\end{figure}

In case of uniform meshes, from equation \eqref{eq_matrix_uniform_grid}, we have that $A_N=c{T}_N(p^{\beta}_N(\theta))$ where $c=\frac{Kh^{\beta-1}}{2^{\beta}\Gamma(\beta+1)}$ and
\begin{equation}\label{eq_symbol}          
p^{\beta}_N(\theta)=\frac{1}{c}a_{1,1}+\frac{2}{c}\sum\limits_{k=1}^{N-1}a_{1,k+1}\cos(k\theta),
\end{equation}
with $a_{1,k}$ in equation \eqref{eq_matrix_uniform_grid}. \credd{The scaling $\frac{1}{c}$ in $p^{\beta}_N(\theta)$ is considered in order to directly deduce the following proposition from Proposition 3.15 in \cite{FVE_RL_marco}}.
\begin{proposition}[\cite{FVE_RL_marco}]\label{prop_spectral_info}
	For $N\rightarrow\infty$, $p^{\beta}_N(\theta)$ in \eqref{eq_symbol} converges to a positive real-valued even function, say $p^{\beta}(\theta)$, that has a unique zero at $\theta=0$ of order lower than $2$, for every $\beta\in(0,1)$ and such that
	$$\{h^{1-\beta}A_N\}_{N}\sim_\lambda \left( \frac{K}{2^{\beta}\Gamma(\beta+1)}p^\beta(\theta), [-\pi,\pi]\right).$$	
\end{proposition}
\credd{In case of non-uniform meshes as in equation \eqref{eq_graded_mesh}, the following theorem, proved in appendix \ref{appendix_1}, holds.}
\begin{theorem}\label{th_symbol_val_sing}
	Let $K(x)=K$ and suppose $g:[0,1]\rightarrow[0,1]$ is an increasing bijective map in $C^3([0, 1])$. Then, if $g'({x})$ has a finite amount of zeros of limited order, it holds
	\begin{equation}\label{eq_temp_th_valsing}
	\{ h^{1-\beta}A_{N}\}_{N} \sim_{\sigma} \left(f_\beta(x,\theta),[0,1]\times[-\pi,\pi]\right),
	\end{equation}
	{with
		\begin{equation}\label{eq_temp_symbol_f}
		f_\beta(x,\theta)=\frac{K}{2^\beta \Gamma(\beta+1) (g'(x))^{1-\beta}}p^{\beta}(\theta).
		\end{equation}
	}
\end{theorem}
Theorem \ref{th_symbol} extends the result obtained in Theorem \ref{th_symbol_val_sing} to the distribution in the sense of the eigenvalues.
\begin{theorem}\label{th_symbol}
	Let $K(x)=K$ and suppose $g:[0,1]\rightarrow[0,1]$ is an increasing bijective map in $C^3([0, 1])$. Then, if $g'({x})>0$, it holds
	\begin{equation}\label{eq:dist_eig}
	\{h^{1-\beta}A_{N}\}_{N} \sim_{\lambda} \left(f_\beta(x,\theta),[0,1]\times[-\pi,\pi]\right),
	\end{equation}
	{with $f_\beta(x,\theta)$ defined as in \eqref{eq_temp_symbol_f}.}
\end{theorem}
\begin{remark}\label{rem_order0}
	Combining Theorem \ref{th_symbol} and Proposition \ref{prop_spectral_info}, we have that the symbol $f_\beta(x,\theta)$ of $\{h^{1-\beta}A_N\}_N$ has a unique zero at $\theta=0$ of order lower than $2$.
\end{remark}
The constraint $g'({x})>0$ is taken to facilitate the proof as, under this hypothesis, we could show that \textbf{GLT5} holds; see appendix \ref{appendix_3} for more details on the proof of Theorem \ref{th_symbol}. Despite the need of this assumption to accomplish the proof, a numerical check indicates that \textbf{GLT5} still holds when $g({x})={x}^q$, as long as 
\begin{equation}\label{eq:q}
1<q<\frac{2-\beta}{1-\beta}. 
\end{equation}
\cred{A confirmation is given in Figure \ref{fig_sequence} where, fixed $K(x)=1$, $\gamma=0.5$, and $g(x)=x^q$, with $q=q_1,q_2$ such that $q_1<\frac{2-\beta}{1-\beta}<q_2$, we plot the sequence
	$$s(N)=\frac{h^{1-\beta}\norm{A_{N}-A_{N}^\mathrm{H}}_\mathrm{tr}}{N},$$	
	varying $N\in[2^4,2^{10}]$ and $\beta$. When $q=q_1$ (resp. $q=q_2$), the sequence $s(N)$ shows a monotonically decreasing (resp. increasing) character, which confirms that \textbf{GLT5} still holds when \eqref{eq:q} is satisfied. A further evidence is given in Figure \ref{fig_GLT5_q} where we 
	plot the sign of $s(2^4)-s(2^5)$ 
	varying both $\beta$ and $q$.
	At each coordinate $(\beta,q)$ the blue dot means that $s(2^4)-s(2^5)>0$ while the yellow dot corresponds to $s(2^4)-s(2^5)<0$.
	If, in line with the results in Figure~\ref{fig_sequence}, we assume that $s(N)$ is monotonic, the blue and yellow areas roughly indicate for which pairs $(\beta,q)$ the sequence $s(N)$ converges or diverges. Note that the curve $\frac{2-\beta}{1-\beta}$, depicted in red, delimits the two regions. When $\beta\approx 0.9$ the red curve does not properly overlap the blue dots and this could be due to the large values of $q$, which lead to highly ill-conditioned linear systems even with small $N$.
}

\begin{figure}
	\centering
	\begin{subfigure}[b]{0.48\textwidth}
		\centering
		\includegraphics[width=\textwidth]{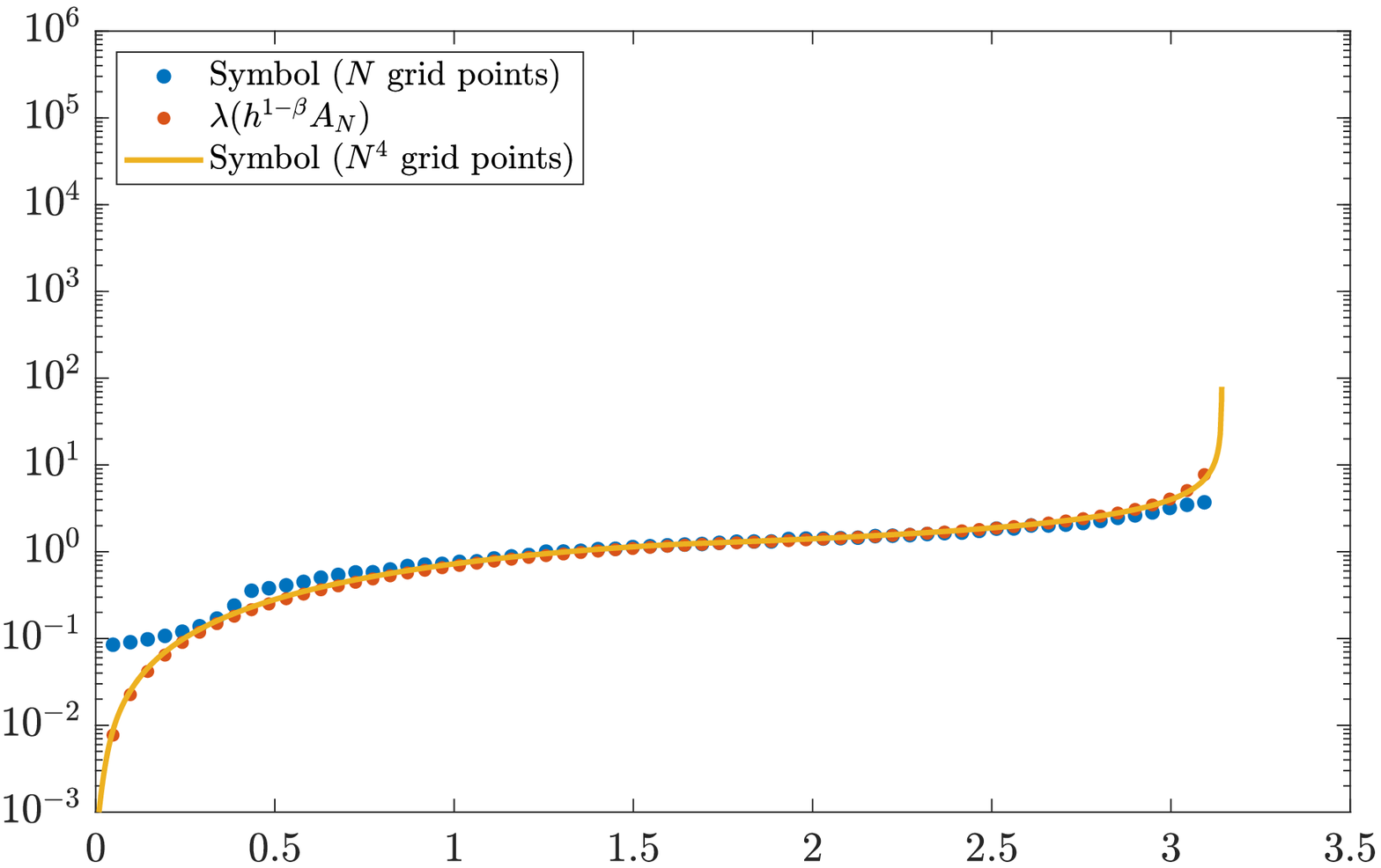}
		\caption{$g(x)=x^2$.}
		\label{fig_non_unif_distribution1}
	\end{subfigure}
	\hfill
	\begin{subfigure}[b]{0.48\textwidth}
		\centering
		\includegraphics[width=\textwidth]{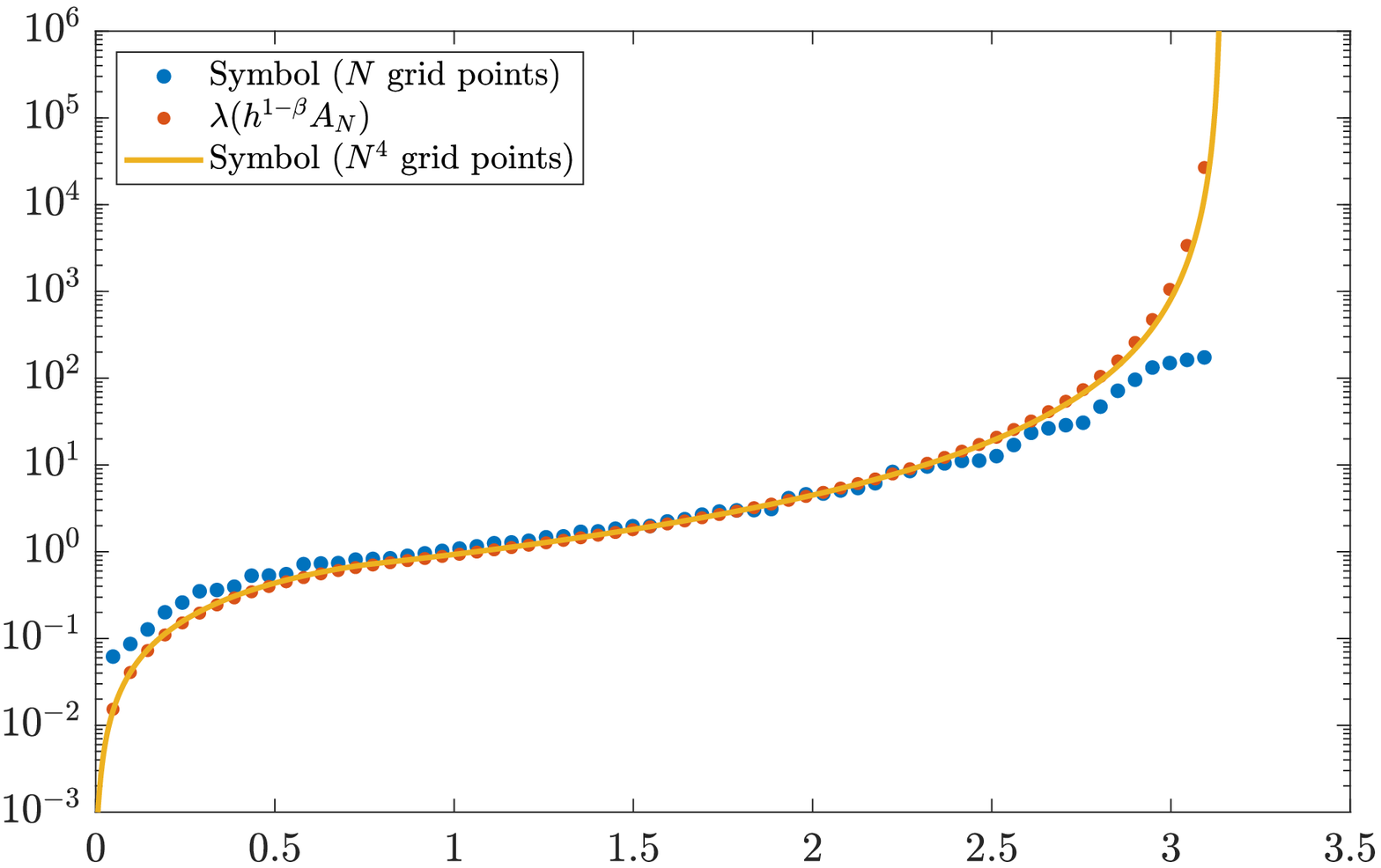}
		\caption{$g(x)=x^4$.}
		\label{fig_non_unif_distribution2}
	\end{subfigure}
	\caption{Plot of the eigenvalues of \cred{$h^{1-\beta}A_N$} with $\beta=0.5$, $N=2^6$ and of a sampling of $f_\beta(x,\theta)$ over the meshes mapped by $g(x)=x^q$, $q=2,4$ from the grids given in (i)--(ii).}\label{fig_non_unif_distribution}
\end{figure}

A further check confirms that Theorem \ref{th_symbol} seems to hold even when \eqref{eq:q} is not satisfied and \textbf{GLT5} does not hold. 
Fixed $N=2^6$ and $\beta=0.5$, Figures \ref{fig_non_unif_distribution1} and \ref{fig_non_unif_distribution2} compare the sorted eigenvalues of $h^{1-\beta}A_N$ with the sorted uniform sampling of the symbol $f_\beta(x,\theta)$ over the meshes mapped by $g({x})={x}^q$, $q=2,4$ from both
\begin{itemize}
	\item[(i)] $\lbrace \hat{x}_i,\theta_j\rbrace_{i,j=1}^{\sqrt{N}}$ with $\hat{x}_i=i\frac{1}{\sqrt{N}}$, $\theta_j=j\frac{\pi}{\sqrt{N}+1}$,
	\item[(ii)] $\lbrace \hat{x}_i, \theta_j\rbrace_{i,j=1}^{N^2}$ with $\hat{x}_i=i\frac{1}{N^2}$, $\theta_j=j\frac{\pi}{N^2+1}$.
\end{itemize}
Note that for $\beta=0.5$, only $q=2$ satisfies condition \eqref{eq:q}. Despite this, in both Figures \ref{fig_non_unif_distribution1} and \ref{fig_non_unif_distribution2} we observe a similar shape between the eigenvalues of \cred{$h^{1-\beta}A_N$} and the sampling of the symbol over the grid in (i), with a lack of overlapping at initial and final grid points. This discrepancy is immediately overcome by making a comparison between the eigenvalues of \cred{$h^{1-\beta}A_N$} and the sampling of the symbol $f_\beta(x,\theta)$ over the much finer mesh in (ii). Further tests, not reported here, show that $f_\beta(x,\theta)$ keeps approximating the eigenvalues distribution of \cred{$h^{1-\beta}A_N$} also for $\beta\neq 0.5$ and $q>\frac{2-\beta}{1-\beta}$. Such a result is in line with Theorem 10.11 in \cite{serra_vol1}, which states that the GLT symbol of equation \eqref{eq_FDE} with $\beta=0$, i.e., the Laplacian, is the spectral symbol independently of the grading parameter of the map $g(x)$.

\section{Multigrid methods} \label{sec:vcycle}
Multigrid methods, introduced in \cite{ruge_stuben} to deal with positive definite linear systems, combine two iterative methods known as smoother and Coarse Grid Correction (CGC). The smoother is typically a simple stationary iterative method. The multigrid algorithm can be figured out starting from the two-grid case. One step of a two-grid method is obtained by: 1) computing an initial approximation by few iterations of a pre-smoother, 2) projecting and solving the error equation into a coarser grid, 3) interpolating the solution of the coarser problem, 4) updating the initial approximation, and finally 5) applying a few iterations of a post-smoother to further improve the approximation. Since the coarser grid could be too large for a direct computation of the solution, the same idea can be recursively applied obtaining the so-called V-cycle method.

A common approach to define the coarser operator, known as \emph{geometric approach}, consists in rediscretizing the same problem on the coarser grid. This approach has the advantage of {preserving} the same structure of the coefficient matrix at each level. On the other hand, the coarser problems need to be properly scaled and the result is usually less robust than the so-called \emph{Galerkin approach}. The latter, for a given linear system $A_{N}x=b$, $A_N\in\mathbb{R}^{N\times N}$, defines the coarser matrix as $A_K=P_N^TA_NP_N$, where $P_N\in\mathbb{R}^{N\times K}$ is the full-rank prolongation matrix, while $P_N^T$ is the restriction operator. The Galerkin approach is useful for the convergence analysis, but in practice it could be computationally too expensive for FDE problems.

The convergence of the V-cycle, {for $A_N$ positive definite}, relies on the so-called \emph{smoothing property} and \emph{approximation property} (see \cite{ruge_stuben}). In order to discuss the convergence analysis of V-cycle applied to \eqref{eq_LINEAR_system_discretizedeq}, we consider the mesh to be uniform, the diffusion coefficient $K$ to be constant, $\gamma=\frac{1}{2}$, and we use weighted Jacobi as smoother.
Under these assumptions and because of the Toeplitz structure of the considered matrices, the weighted Jacobi is well-known to satisfy the smoothing property for positive definite matrices, whenever it is convergent \cite{smoothing_p_cite}. {Moreover, according to Proposition \ref{prop_spectral_info}, in case of uniform meshes the symbol \cred{$f_\beta(x,\theta)$ of $h^{1-\beta}A_N$} vanishes with order lower than $2$ at $\theta=0$ and the approximation property holds true with the same projectors as in the case of the Laplacian (see \cite{M2D2}). However, in case of non-uniform meshes, the grid transfer operator should take into account the nonuniformity of the grid, see \cite{briggs}.
	
	
	We propose a V-cycle whose hierarchy is built through the geometric approach. Assuming that $\Omega_\ell=\lbrace x_i^\ell\rbrace_{i=0}^{N^\ell+1}$ is the grid on the $\ell$-th level, where $\ell\geq0$ and $\Omega_0$ is the finest mesh, then $\Omega_{\ell+1},\ \ell\geq0,$ is obtained from $\Omega_{\ell}$ as
	\begin{equation*}
	\Omega_{\ell+1}=
	\{x_0^{\ell},x_{N^\ell+1}^{\ell}\}\cup\{x_{2k}^{\ell}\}_{k=1}^{N^{\ell+1}},
	\end{equation*}
	where $N^{\ell+1}=\lfloor \frac{N^\ell}{2}\rfloor$ are the degrees of freedom of equation \eqref{eq_FDE} discretized over $\Omega_{\ell+1}$. By iterating this process $\mathrm{lvl}$ times, with $\Omega_{\mathrm{lvl}}$ such that $N^\mathrm{lvl}\leq3$, we obtain the V-cycle hierarchy. Note that, in order to make the V-cycle properly working, the linear systems must be scaled such that the right-hand side does not contain any grid dependent scaling factor. Therefore, we multiply both members of $A_Nx=b$ by the diagonal matrix $H_N=\diag\limits_{i=1,...,N}(\frac{1}{h_i})$. Finally, the projectors are built according to the discussion in \cred{Chapter 7 at page 129 of \cite{briggs}.}

	Regarding the smoother, at each iteration of V-cycle one iteration of relaxed Jacobi as pre- and post-smoother is performed. The relaxation parameter $\omega$ is estimated through the approach introduced in \cite{paper_MGS}. Such estimation is obtained by: 1) rediscretizing equation \eqref{eq_FDE} over a coarser grid ($\tilde{N}\!\leq\! 2^4$), 2) computing the spectrum of the Jacobi iteration matrix, 3) choosing the weight $\omega$ in such a way that the whole spectrum is contained inside a complex set $O=\lbrace (x,y)\ |\ x\in I\subset\mathbb{R}, -\tilde{o}(x)<y<\tilde{o}(x)\rbrace$. A possible choice for $\tilde{o}(x)$ is given by $\tilde{o}(x)=\sqrt{1-x^2}+c x-c,\ c>0$, which is the sum of a semicircle and a line, and is motivated by the need of clustering the spectrum of the Jacobi iteration matrix inside the unitary circle. Note that $\tilde{o}(x)$ yields a set $O$ that is slightly smaller than the unitary circle in such a way that possible outliers are still smaller than 1 in modulus. Our numerical tests in Section \ref{sec_results} confirm the suitability of the choice $\tilde{o}(x)=\sqrt{1-x^2}+0.475x-0.475$ and $I=\left[-\tfrac{1239}{1961},1\right]\approx[-0.63,1]$ {such that $\tilde{o}(x)\geq 0$}. Note that our choice slightly differs from the one proposed in \cite{paper_MGS}, where the authors dealt with a different model.
	
	
	
	\begin{remark}\label{rem:laplacian}
		From Remark \ref{rem_b0_b1}, in the case where $\beta=0$, $A_N$ becomes the positive definite one-dimensional discrete Laplacian, therefore we expect that our multigrid performs well for $\beta\approx 0$ even in the anisotropic cases where $\gamma\approx0$ or $\gamma\approx1$. When $\beta=1$, $A_N$ becomes skew-symmetric, hence we do not expect that our multigrid method performs well in the cases where $\beta\approx 1$ and $\gamma\approx0$ or $\gamma\approx1$.
	\end{remark}

	\section{Numerical results}\label{sec_results}
	In this section, we compare few (graded and composite) grids by reporting the reconstruction error and the convergence order and we check the performances of our multigrid applied to equation \eqref{eq_FDE} varying the grid, $\gamma$ and $\beta$. Precisely, aiming at increasing its robustness, we use the multigrid method described in Section \ref{sec:vcycle} as preconditioner for GMRES by performing one iteration of V-cycle applied directly to the coefficient matrix. Throughout, we denote our solver by $\mathcal{P}$-GMRES.
	
	Our numerical tests have been run on a server with AMD 3600 6-core (4.20 GHz) processor and 64 GB (3600 MHz) RAM and Matlab 2020b. In all our tests we consider equation \eqref{eq_FDE} with $f(x)=\frac{(1-\gamma)(1-\beta)}{\Gamma(\beta)x(1-x)^{1-\beta}}$, $K(x)=1$ with $\gamma\in[0,1]$ and $\beta\in(0,1)$, whose exact solution, according to \cite{wang}, is $u(x)=x^{1-\beta}$.

	For all involved iterative methods the initial guess $x^{(0)}$ is the null vector, the maximum amount of iterations is $100$ and the stopping criterion is 
	\begin{equation*}
	\frac{\norm{Ax^{(k)}-b}_2}{\norm{b}_2}<\mathrm{tol},
	\end{equation*}
	where the tolerance is tol$=10^{-7}$ and $x^{(k)}$ is the unknown at the $k$-th iteration.\\

	\paragraph{Test 1} In this first test, we discuss the choice of the parameter $q$ for the function $g_{q,\mathbf{\epsilon}}(x)$ in \eqref{eq_grid_g} which generates the {graded meshes} discussed in Section \ref{sec_grid}. Precisely, fixed $N=2^{10}$, we compare in terms of approximation error the value $q_\beta$, defined in \eqref{eq_q} and given in \cite{Kopteva}, with the numerically computed optimal value $q_{\text{opt}}$ obtained by minimizing the infinity norm error on a set of equispaced values in the interval $[1,9]$. 
	
	Table \ref{tab_q_vary_grid} shows the optimal value $q_{\text{opt}}$ and the infinity norm errors $e^{\text{opt}}$ and $e^\beta$ yield when discretizing equation \eqref{eq_FDE} over the grid mapped by $g_{q,\mathbf{\epsilon}}(x)$ with $q=q_{\text{opt}}$ and $q=q_\beta$, respectively, for various choices of $\mathbf{\epsilon}=(\epsilon_1,\epsilon_2)$. Note that $q_{0.2}=1.5,\ q_{0.5}=3$ are obtained from \eqref{eq:q}, while the case $\beta=0.8$ deserves a special attention as $q_{0.8}=\frac{1+0.8}{1-0.8}=9$ yields a step length much {lower than} the machine precision for $N=2^{10}$ ($h_1\approx 10^{-28}$). Therefore, according to the discussion in Section \ref{sec_grid}, we set $q_{0.8}=5.3$ such that $h_1= 10^{-16}$.
	
	When $\beta=0.2$, in Table \ref{tab_q_vary_grid} we observe an increase in the error when choosing $q_\beta$ with respect to $q_{\text{opt}}$, i.e., $e^\beta\approx 2e^{\text{opt}}$ independently of the choice of $\gamma$ and of $\mathbf{\epsilon}$. Moreover, the interval lengths $\epsilon_1,\epsilon_2$ do not seem to affect the error in the case of $q_\beta$ (all $\mathbf{\epsilon}^{(i)}$ yield the same $e^\beta$). This also means that the use of a smooth function (cases $\mathbf{\epsilon}^{(1)}$, $\mathbf{\epsilon}^{(2)}$, and $\mathbf{\epsilon}^{(4)}$) for generating the grid does not decrease the error in comparison with a non-smooth function (cases $\mathbf{\epsilon}^{(3)}$, $\mathbf{\epsilon}^{(5)}$, and $\mathbf{\epsilon}^{(6)}$), therefore smaller values of $\epsilon_1$ and $\epsilon_2$ can be chosen, in order to speed up the matrix-vector product. When $\beta=0.5$, the difference between $e^\beta$ and $e^{\text{opt}}$ is almost negligible, both do not vary that much with $\gamma$ and $\epsilon$, and again the choice of $\epsilon_1,\epsilon_2$ is not so crucial. When $\beta=0.8$, unlike the previous two cases, if the non linear part of function $g_{q,\mathbf{\epsilon}}(x)$ is too short, e.g., when considering $\mathbf{\epsilon}^{(1)}$, the error increases greatly. This happens because more grid points are needed near $x=0$ to deal with the {singularity} of the solution, and by increasing $\epsilon_1$, $g_{q,\mathbf{\epsilon}}(x)$ projects more grid points near $x=0$.
	
	Moreover, by comparing $\mathbf{\epsilon}^{(2)}$ with $\mathbf{\epsilon}^{(3)}$ we note that when $g_{q,\mathbf{\epsilon}}(x)$ is smooth, i.e., $\mathbf{\epsilon}=\mathbf{\epsilon}^{(2)}$, the error is lower than in the non-smooth case $\mathbf{\epsilon}=\mathbf{\epsilon}^{(3)}$, especially when $\beta\geq 0.5$, even if the length $\epsilon_1+\epsilon_2$ of the interval where $g_{q,\mathbf{\epsilon}}(x)$ is non linear does not change.  Therefore, a smooth function $g_{q,\mathbf{\epsilon}}(x)$ is recommended when $\beta\approx 1$.
	
	Note that, in any of the tested cases, a full non-uniform mesh does not seem to be necessary, since the lowest error is always reached for mixed meshes.
	
	\begin{table}[!ht]\scriptsize\begin{center}
			\begin{tabular}{@{\extracolsep{3pt}}c@{\hspace{0.2cm}}c@{\hspace{0.2cm}}c@{\hspace{0.2cm}}c@{\hspace{0.1cm}}l@{\hspace{0.1cm}}c@{\hspace{0.15cm}}c@{\hspace{0.1cm}}c@{\hspace{0.1cm}}c@{\hspace{0.15cm}}c@{\hspace{0.1cm}}c@{\hspace{0.1cm}}c@{\hspace{0.15cm}}c@{\hspace{0.1cm}}c@{\hspace{0.1cm}}c@{\hspace{0.15cm}}c@{\hspace{0.1cm}}c@{\hspace{0.1cm}}c@{\hspace{0.15cm}}c@{\hspace{0.1cm}}c@{\hspace{0.1cm}}c@{\hspace{0.15cm}}}
				\hline\vspace{-0.2cm} \\
				\multirow{2}{*}{$\gamma$} &\multirow{2}{*}{$\beta$ } &\multirow{2}{*}{($q_\beta$)}  & \multicolumn{3}{c}{$\mathbf{\epsilon}^{(1)}=(0.1,0.05)$} & \multicolumn{3}{c}{$\mathbf{\epsilon}^{(2)}=(0.2,0.05)$} & \multicolumn{3}{c}{$\mathbf{\epsilon}^{(3)}=(0.25,0)$} & \multicolumn{3}{c}{$\mathbf{\epsilon}^{(4)}=(0.45,0.05)$} & \multicolumn{3}{c}{$\mathbf{\epsilon}^{(5)}=(0.5,0)$} & \multicolumn{3}{c}{$\mathbf{\epsilon}^{(6)}=(1,0)$}\\\cline{4-6}\cline{7-9}\cline{10-12}\cline{13-15}\cline{16-18}\cline{19-21}
				\vspace{-0.2cm} \\
				& & & $q_\text{opt}$ & $e^\text{opt}$ & $e^\beta$ & $q_\text{opt}$ & $e^\text{opt}$ & $e^\beta$ & $q_\text{opt}$ & $e^\text{opt}$ & $e^\beta$ & $q_\text{opt}$ & $e^\text{opt}$ & $e^\beta$ & $q_\text{opt}$ & $e^\text{opt}$ & $e^\beta$ & $q_\text{opt}$ & $e^\text{opt}$ & $e^\beta$ \\ \hline\hline
				\vspace{-0.2cm} \\
				\multirow{3}{*}{$0.3$}
				& $0.2$ & $(1.5)$  & {1.9}& {1.2e-5}& {3.3e-5} & {1.8}& {1.3e-5}& {3.3e-5} & {1.8}& {1.3e-5}& {3.3e-5} & {1.7}& {1.5e-5}& {3.3e-5} & {1.7}& {1.5e-5}& {3.3e-5} & {1.7}& {1.8e-5}& {3.3e-5} \\
				& $0.5$ &$(3.0)$  & {2.7}& {5.2e-5}& {5.7e-5} & {3.0}& {2.6e-5}& {2.6e-5} & {2.9}& {3.1e-5}& {3.2e-5} & {3.1}& {2.1e-5}& {2.2e-5} & {3.0}& {2.2e-5}& {2.2e-5} & {2.8}& {3.5e-5}& {3.7e-5} \\
				& $0.8$ &(5.3)  & {4.4}& {1.9e-3}& {6.3e-3} & {5.3}& {7.1e-4}& {7.1e-4} & {4.6}& {1.6e-3}& {5.7e-3} & {5.3}& {7.1e-4}& {7.1e-4} & {5.3}& {7.1e-4}& {7.1e-4} & {5.3}& {7.1e-4}& {7.1e-4} \\
				\hline\vspace{-0.2cm} \\  
				\multirow{3}{*}{$0.5$}  
				& $0.2$ & (1.5) & {1.7}& {9.0e-6}& {2.1e-5} & {1.7}& {9.6e-6}& {2.1e-5} & {1.7}& {9.7e-6}& {2.1e-5} & {1.7}& {1.1e-5}& {2.1e-5} & {1.7}& {1.1e-5}& {2.1e-5} & {1.6}& {1.3e-5}& {2.1e-5} \\
				& $0.5$ & (3.0) & {2.9}& {1.3e-5}& {1.4e-5} & {3.0}& {1.1e-5}& {1.1e-5} & {3.0}& {1.2e-5}& {1.2e-5} & {2.9}& {1.5e-5}& {1.5e-5} & {2.9}& {1.5e-5}& {1.5e-5} & {2.7}& {2.4e-5}& {2.6e-5} \\
				& $0.8$ & (5.3) & {4.3}& {1.0e-3}& {4.2e-3} & {5.3}& {3.5e-4}& {3.7e-4} & {4.6}& {8.0e-4}& {2.9e-3} & {5.3}& {3.4e-4}& {3.4e-4} & {5.3}& {3.4e-4}& {3.4e-4} & {5.3}& {3.4e-4}& {3.4e-4} \\
				\hline\vspace{-0.2cm} \\  
				\multirow{3}{*}{$0.7$}  
				& $0.2$& (1.5) & {1.7}& {5.6e-6}& {1.3e-5} & {1.7}& {6.0e-6}& {1.3e-5} & {1.7}& {6.1e-6}& {1.3e-5} & {1.7}& {6.9e-6}& {1.3e-5} & {1.7}& {7.0e-6}& {1.3e-5} & {1.6}& {8.1e-6}& {1.3e-5} \\
				& $0.5$& (3.0) & {2.6}& {1.4e-5}& {1.8e-5} & {2.7}& {8.2e-6}& {1.0e-5} & {2.7}& {9.8e-6}& {1.5e-5} & {2.7}& {8.5e-6}& {8.7e-6} & {2.7}& {8.7e-6}& {9.0e-6} & {2.6}& {1.3e-5}& {1.5e-5} \\
				& $0.8$& (5.3) & {3.8}& {2.1e-4}& {1.9e-3} & {4.1}& {4.4e-5}& {4.6e-4} & {3.7}& {2.5e-4}& {5.3e-3} & {4.2}& {3.7e-5}& {5.1e-4} & {4.2}& {3.7e-5}& {5.5e-4} & {4.2}& {3.7e-5}& {5.2e-4} \\
				\hline
			\end{tabular}
		\end{center}
		\caption{Comparison between the numerically computed optimal value $q_{\text{opt}}$ and $q_\beta$, in terms of infinity norm approximation error varying $\gamma,\beta$ and~$\mathbf{\epsilon}$.}\label{tab_q_vary_grid}
	\end{table}

	\paragraph{Test 2} 
	We now fix $\gamma=0.5$ and test the robustness of $\mathcal{P}$-GMRES by reporting the iterations to tolerance (It), needed for solving equation \eqref{eq_FDE} discretized over the grids reported in Section \ref{sec_grid}. 
	
	Table \ref{tab_ex2} shows It, as well as the numerical infinity norm error $e^\infty$ and the convergence order computed as the ratio between the infinity-norm error over two grids, with $N$ and $2N$ points, in $\log_2$ scale. With symbol \lq\lq -'' we mean that the solver exceeded the maximum amount of iterations fixed to $100$. 
	
	We observe that, when considering a graded mesh mapped by $g_{q,\mathbf{\epsilon}}(x)$, for all choices of $\beta$, $\mathcal{P}$-GMRES converges linearly. On the other hand, for both given composite meshes and for $\beta\le0.5$, It increases with $N$ and the linear convergence is lost. {Less} iterations are obtained for $\beta=0.8$, but the approximation error yield by the composite meshes substantially worsen compared to the one of the meshes mapped by $g_{q,\mathbf{\epsilon}}(x)$.
		
		In terms of accuracy, when considering $\beta=0.2$, the convergence order of the grids mapped by $g_{q,\mathbf{\epsilon}}$ is $1+\beta$, which coincides with the theoretical convergence rate obtained in \cite{Kopteva}. When $\beta=0.5$, the mapping function $g_{q,\mathbf{\epsilon}}(x)$ with $\mathbf{\epsilon}=\mathbf{\epsilon}^{(1)}$ seems to have a too short non linear part, since for ${N<2^{9}-1}$ it yields a larger error in comparison with other choices of $\mathbf{\epsilon}$. {When $\beta=0.8$, for
		\begin{itemize}
			\item $\mathbf{\epsilon}=\mathbf{\epsilon}^{(6)}$ the order $\beta+1=1.8$ is obtained only for $N\leq 2^6-1$, while for larger $N$ the order decreases. This is caused by the {lower cap \eqref{eq_lower_cap}} imposed on the smallest step size, i.e., $h_1=10^{-16}$, which is needed to avoid stability problems related to the machine precision;
			\item $\mathbf{\epsilon}=\mathbf{\epsilon}^{(4)}$, similarly to $\mathbf{\epsilon}=\mathbf{\epsilon}^{(6)}$, the convergence order decreases with $N$ for $N>2^8-1$. When $2^7-1\leq N\leq2^8-1$ we have $ord>1+\beta$, but $e^\infty$ for $2^6-1\leq N\leq2^7-1$ is still much larger than in the case of $\mathbf{\epsilon}^{(6)}$, therefore $g_{q,\mathbf{\epsilon}}$ with $\mathbf{\epsilon}=\mathbf{\epsilon}^{(4)}$ does not allow higher convergence order \credd{than with $\mathbf{\epsilon}=\mathbf{\epsilon}^{(6)}$};
			\item $\mathbf{\epsilon}=\mathbf{\epsilon}^{(2)}$, which yields a mesh with a shorter non-uniform part compared the one yield by $\mathbf{\epsilon}=\mathbf{\epsilon}^{(4)}$, the convergence order keeps increasing and at $N=2^{10}-1$ \credd{provides} almost the same accuracy error than $g_{q,\mathbf{\epsilon}}$ with $\mathbf{\epsilon}=\mathbf{\epsilon}^{(4)},\mathbf{\epsilon}^{(6)}$.
		\end{itemize}
	}
		This suggests that the length of the non linear part, i.e., $\epsilon_1+\epsilon_2$, should depend not only on $\beta$, as already observed in Test 1, but also on $N$: the larger $N$ is, the shorter $\epsilon_1+\epsilon_2$ can be. 
	
	
	
	\begin{table}[!ht]\scriptsize\begin{center}
			\begin{tabular}{@{\extracolsep{3pt}}c@{\hspace{0.2cm}}c@{\hspace{0.2cm}}g@{\hspace{0.2cm}}c@{\hspace{0.1cm}}c@{\hspace{0.2cm}}g@{\hspace{0.2cm}}c@{\hspace{0.1cm}}c@{\hspace{0.2cm}}g@{\hspace{0.2cm}}c@{\hspace{0.1cm}}c@{\hspace{0.2cm}}g@{\hspace{0.2cm}}c@{\hspace{0.1cm}}c@{\hspace{0.2cm}}g@{\hspace{0.2cm}}c@{\hspace{0.1cm}}c@{\hspace{0.2cm}}g@{\hspace{0.2cm}}c@{\hspace{0.1cm}}c@{\hspace{0.2cm}}g@{\hspace{0.2cm}}c@{\hspace{0.1cm}}c@{\hspace{0.1cm}}c@{\hspace{0.15cm}}c@{\hspace{0.1cm}}}
				\hline\vspace{-0.2cm} \\
				\multirow{4}{*}{$\beta$} &\multirow{4}{*}{$N\!+\!1$}& \multicolumn{6}{c}{Composite mesh mapped by $g_G(N)$}&\multicolumn{12}{c}{Non-uniform mesh mapped by $g_{q,\mathbf{\epsilon}}(x)$}\\\cline{3-8}\cline{9-20}\vspace{-0.2cm}\\
				&  & \multicolumn{3}{c}{$g_G(N)=\lfloor \sqrt{N}\rfloor$} & \multicolumn{3}{c}{$g_G(N)=\lfloor \log_2{N}\rfloor$} & \multicolumn{3}{c}{$\mathbf{\epsilon}^{(1)}=(0.1,0.05)$} & \multicolumn{3}{c}{$\mathbf{\epsilon}^{(2)}=(0.2,0.05)$} & \multicolumn{3}{c}{$\mathbf{\epsilon}^{(4)}=(0.45,0.05)$} & \multicolumn{3}{c}{$\mathbf{\epsilon}^{(6)}=(1,0)$}\\\cline{3-5}\cline{6-8}\cline{9-11}\cline{12-14}\cline{15-17}\cline{18-20}
				\vspace{-0.2cm} \\
				& & It & $e^\infty$ & ord & It & $e^\infty$ & ord & It & $e^\infty$ & ord & It & $e^\infty$ & ord & It & $e^\infty$ & ord & It & $e^\infty$ & ord \\ \hline\hline
				\multirow{7}{*}{$0.2$}
				& $2^{4}$  & {11}& {2.9e-3}& {} & {11}& {2.9e-3}& {} & {8}& {3.3e-3}& {} & {10}& {3.0e-3}& {} & {10}& {3.0e-3}& {} & {10}& {3.0e-3}& {} \\
				& $2^{5}$  & {12}& {1.0e-3}& {1.5} & {11}& {1.2e-3}& {1.2} & {8}& {1.4e-3}& {1.3} & {8}& {1.3e-3}& {1.2} & {8}& {1.3e-3}& {1.1} & {10}& {1.3e-3}& {1.1} \\
				& $2^{6}$  & {16}& {4.7e-4}& {1.1} & {11}& {5.8e-4}& {1.1} & {8}& {5.9e-4}& {1.2} & {8}& {5.9e-4}& {1.2} & {8}& {5.9e-4}& {1.2} & {8}& {5.9e-4}& {1.2} \\
				& $2^{7}$  & {27}& {2.5e-4}& {0.9} & {17}& {2.9e-4}& {1.0} & {8}& {2.6e-4}& {1.2} & {8}& {2.6e-4}& {1.2} & {8}& {2.6e-4}& {1.2} & {8}& {2.6e-4}& {1.2} \\
				& $2^{8}$  & {35}& {1.4e-4}& {0.8} & {22}& {1.5e-4}& {0.9} & {8}& {1.1e-4}& {1.2} & {8}& {1.1e-4}& {1.2} & {8}& {1.1e-4}& {1.2} & {8}& {1.1e-4}& {1.2} \\
				& $2^{9}$  & {40}& {8.1e-5}& {0.8} & {28}& {8.4e-5}& {0.9} & {8}& {4.9e-5}& {1.2} & {8}& {4.9e-5}& {1.2} & {8}& {4.9e-5}& {1.2} & {9}& {4.9e-5}& {1.2} \\
				& $2^{10}$  & {-}& {-}& {-} & {36}& {4.7e-5}& {0.8} & {8}& {2.1e-5}& {1.2} & {8}& {2.1e-5}& {1.2} & {8}& {2.1e-5}& {1.2} & {9}& {2.1e-5}& {1.2} \\
				\hline\vspace{-0.2cm} \\  
				\multirow{7}{*}{$0.5$}  
				& $2^{4}$  & {8}& {2.3e-2}& {} & {8}& {2.3e-2}& {} & {7}& {2.0e-2}& {} & {8}& {1.2e-2}& {} & {10}& {5.1e-3}& {} & {10}& {5.0e-3}& {} \\
				& $2^{5}$  & {9}& {8.3e-3}& {1.5} & {9}& {1.1e-2}& {1.0} & {8}& {9.4e-3}& {1.1} & {8}& {4.1e-3}& {1.5} & {9}& {1.8e-3}& {1.5} & {10}& {1.8e-3}& {1.5} \\
				& $2^{6}$  & {11}& {2.9e-3}& {1.5} & {11}& {5.7e-3}& {1.0} & {9}& {3.2e-3}& {1.6} & {9}& {1.2e-3}& {1.8} & {9}& {6.4e-4}& {1.5} & {10}& {6.4e-4}& {1.5} \\
				& $2^{7}$  & {13}& {1.3e-3}& {1.1} & {11}& {2.8e-3}& {1.0} & {9}& {9.0e-4}& {1.8} & {9}& {3.1e-4}& {1.9} & {9}& {2.3e-4}& {1.5} & {10}& {2.3e-4}& {1.5} \\
				& $2^{8}$  & {17}& {9.1e-4}& {0.6} & {13}& {1.4e-3}& {1.0} & {10}& {2.3e-4}& {2.0} & {10}& {8.0e-5}& {1.9} & {10}& {8.0e-5}& {1.5} & {10}& {1.0e-4}& {1.1} \\
				& $2^{9}$  & {22}& {6.4e-4}& {0.5} & {16}& {7.8e-4}& {0.9} & {10}& {5.6e-5}& {2.0} & {11}& {2.8e-5}& {1.5} & {11}& {3.0e-5}& {1.4} & {10}& {5.2e-5}& {1.0} \\
				& $2^{10}$  & {22}& {4.5e-4}& {0.5} & {17}& {5.1e-4}& {0.6} & {11}& {1.4e-5}& {2.0} & {10}& {1.1e-5}& {1.3} & {10}& {1.5e-5}& {1.0} & {11}& {2.6e-5}& {1.0} \\
				\hline\vspace{-0.2cm} \\  
				\multirow{7}{*}{$0.8$}
				& $2^{4}$  & {8}& {1.1e-1}& {} & {8}& {1.1e-1}& {} & {7}& {1.2e-1}& {} & {7}& {1.3e-1}& {} & {7}& {9.0e-2}& {} & {8}& {2.2e-2}& {} \\
				& $2^{5}$  & {8}& {7.6e-2}& {0.6} & {8}& {8.7e-2}& {0.4} & {7}& {1.1e-1}& {0.2} & {7}& {9.7e-2}& {0.4} & {8}& {5.0e-2}& {0.8} & {8}& {6.5e-3}& {1.8} \\
				& $2^{6}$  & {10}& {5.0e-2}& {0.6} & {8}& {6.6e-2}& {0.4} & {7}& {7.7e-2}& {0.5} & {7}& {5.8e-2}& {0.7} & {8}& {2.4e-2}& {1.0} & {8}& {1.9e-3}& {1.8} \\
				& $2^{7}$  & {10}& {2.5e-2}& {1.0} & {9}& {5.0e-2}& {0.4} & {7}& {4.5e-2}& {0.8} & {7}& {3.5e-2}& {0.7} & {8}& {4.0e-3}& {2.6} & {8}& {9.2e-4}& {1.1} \\
				& $2^{8}$  & {11}& {1.3e-2}& {1.0} & {9}& {3.8e-2}& {0.4} & {7}& {3.2e-2}& {0.5} & {7}& {1.3e-2}& {1.4} & {8}& {5.8e-4}& {2.8} & {9}& {5.8e-4}& {0.7} \\
				& $2^{9}$  & {13}& {5.5e-3}& {1.2} & {10}& {2.9e-2}& {0.4} & {7}& {1.7e-2}& {0.9} & {8}& {2.7e-3}& {2.3} & {8}& {4.5e-4}& {0.4} & {9}& {4.3e-4}& {0.4} \\
				& $2^{10}$  & {14}& {4.8e-3}& {0.2} & {10}& {2.2e-2}& {0.4} & {7}& {4.2e-3}& {2.0} & {8}& {3.7e-4}& {2.8} & {8}& {3.4e-4}& {0.4} & {9}& {3.4e-4}& {0.3} \\
				\hline
			\end{tabular}
		\end{center}
		\caption{Iterations to tolerance of $\mathcal{P}$-GMRES, approximation error, and convergence order varying $N$ and the grid for $\gamma=0.5$.}\label{tab_ex2}
	\end{table}

	\paragraph{Test 3} We now fix $\gamma=0.5$, $\beta=0.9$ and compare $\mathcal{P}$-GMRES with the Preconditioned Fast Conjugate Gradient Squared (PFCGS) introduced in \cite{wang}, which consists in a T. Chan's block circulant preconditioner. \credd{Precisely, such preconditioner has a two-block structure as a consequence} of the composite mesh combining two type of meshes, the refinement near $x=0$ and the uniform mesh.
	
	Table \ref{tab_pfcgs} shows $e^\infty$ and It of both PFCGS and $\mathcal{P}$-GMRES in case of the composite mesh given in \cite{wang} and described in Section \ref{sec_grid}. We recall that, according to equation \eqref{eq_composite_mesh_x}, $N_1=g_G(N)$ is the number of points of the refined part of the mesh, while $N_2$ is the number of points that compose the uniform part of the mesh and therefore the composite mesh has $N=N_1+N_2$ points. We note that $\mathcal{P}$-GMRES has more stable iteration number when increasing $N$ with respect to PFCGS. 
	
	In Table \ref{tab_ex3}, we only consider $\mathcal{P}$-GMRES and further test it over both composite and graded meshes. Aside from It and $e^\infty$, we also check the 2-norm relative numerical error $e^{\text{rel}}$ of $\mathcal{P}$-GMRES varying $N$ and the grid. Note that, due to the solution being singular, $e^\infty$ describes the error near the singularity, while $e^\mathrm{rel}$ gives an average error on the whole domain.
	
	We note that $\mathbf{\epsilon}^{(6)}$ yields the lowest $e^\infty$ for $2^4\leq N+1\leq 2^6$, but then, {increasing $N$}, $e^\infty$ stops decreasing 
	\crep{while $e^\mathrm{rel}$ still decreases. 
		This is because the {lower cap \eqref{eq_lower_cap}} on the step size 
		does not allow to increase the accuracy near the singularity, where a larger amount of grid points is required, and hence does not allow $e^\infty$ to further decrease. On the other hand, when increasing $N$ the non singular part of the solution is better approximated and hence $e^\mathrm{rel}$ decreases.} Our choice for the lower cap as $h_1=10^{-16}$ could of course be modified. 
	Indeed, we observe that again for ${\bf \epsilon}=\mathbf{\epsilon}^{(6)}$ when $N=2^6-1$ we obtain an $e^\infty$ close to the $e^\infty$ given by the composite mesh with $g_G(N)=\lfloor \sqrt{N}\rfloor$ and $N=2^{10}-1$. Therefore, by improving the lower cap on the step size, the mesh mapped by $g_{q,\mathbf{\epsilon}}(x)$ could potentially lead to lower errors with much smaller sizes compared to the composite mesh.
	
	
	Finally, we note that the iterations are stable as $N$ increases for any of the tested grids, which makes $\mathcal{P}$-GMRES a suitable solver.

	\begin{table}[!ht]\scriptsize\begin{center}
			\renewcommand{\arraystretch}{1.5}
			\begin{tabular}{c@{\hspace{0.5cm}}c@{\hspace{0.5cm}}c@{\hspace{0.5cm}}c@{\hspace{0.5cm}}c}\hline
				$N_1$ & $N_2$ & $e^\infty$ & It PFCGS & It $\mathcal{P}$-GMRES\\
				\hline
				$2^3$ & $2^8$ & $7.9306$e-2 & 13 & 7\\
				$2^4$ & $2^9$ & $4.2326$e-2 & 16 & 7\\
				$2^5$ & $2^{10}$ & $1.3025$e-2 & 25 & 8\\
				\hline
			\end{tabular}
		\end{center}
		\caption{Iterations to tolerance of PFCGS and $\mathcal{P}$-GMRES with $\gamma=0.5$ and $\beta=0.9$.}\label{tab_pfcgs}
	\end{table}

	\begin{table}[!ht]\scriptsize\begin{center}
			\begin{tabular}{@{\extracolsep{3pt}}c@{\hspace{0.2cm}}g@{\hspace{0.2cm}}c@{\hspace{0.1cm}}c@{\hspace{0.2cm}}g@{\hspace{0.2cm}}c@{\hspace{0.1cm}}c@{\hspace{0.2cm}}g@{\hspace{0.2cm}}c@{\hspace{0.1cm}}c@{\hspace{0.2cm}}g@{\hspace{0.2cm}}c@{\hspace{0.1cm}}c@{\hspace{0.2cm}}g@{\hspace{0.2cm}}c@{\hspace{0.1cm}}c@{\hspace{0.2cm}}g@{\hspace{0.2cm}}c@{\hspace{0.1cm}}c@{\hspace{0.2cm}}g@{\hspace{0.2cm}}c@{\hspace{0.1cm}}c@{\hspace{0.1cm}}c@{\hspace{0.15cm}}c@{\hspace{0.1cm}}}
				\hline\vspace{-0.2cm} \\
				& \multicolumn{6}{c}{Composite mesh mapped by $g_G(N)$}&\multicolumn{12}{c}{Non-uniform mesh mapped by $g_{q,\mathbf{\epsilon}}(x)$}\\\cline{2-7}\cline{8-19}\vspace{-0.2cm} \\
				\multirow{2}{*}{$N\!+\!1$}  & \multicolumn{3}{c}{$g_G(N)=\lfloor \sqrt{N}\rfloor$} & \multicolumn{3}{c}{$g_G(N)=\lfloor \log_2{N}\rfloor$} & \multicolumn{3}{c}{$\mathbf{\epsilon}^{(1)}=(0.1,0.05)$} & \multicolumn{3}{c}{$\mathbf{\epsilon}^{(2)}=(0.2,0.05)$} & \multicolumn{3}{c}{$\mathbf{\epsilon}^{(4)}=(0.45,0.05)$} & \multicolumn{3}{c}{$\mathbf{\epsilon}^{(6)}=(1,0)$}\\\cline{2-4}\cline{5-7}\cline{8-10}\cline{11-13}\cline{14-16}\cline{17-19}
				\vspace{-0.2cm} \\
				& It & $e^\infty$ & $e^{\text{rel}}$ & It & $e^\infty$ & $e^{\text{rel}}$ & It & $e^\infty$ & $e^{\text{rel}}$ & It & $e^\infty$ & $e^{\text{rel}}$ & It & $e^\infty$ & $e^{\text{rel}}$ & It & $e^\infty$ & $e^{\text{rel}}$ \\ \hline\hline
				$2^{4}$  & {7}& {1.9e-1}& {8.7e-2} & {7}& {1.9e-1}& {8.7e-2} & {7}& {2.3e-1}& {8.8e-2} & {6}& {2.0e-1}& {9.6e-2} & {6}& {1.8e-1}& {1.2e-1} & {7}& {6.4e-2}& {1.0e-1} \\
				$2^{5}$  & {7}& {1.5e-1}& {5.3e-2} & {7}& {1.6e-1}& {5.6e-2} & {7}& {1.9e-1}& {6.1e-2} & {7}& {1.9e-1}& {6.7e-2} & {7}& {1.2e-1}& {6.4e-2} & {8}& {3.2e-2}& {3.3e-2} \\
				$2^{6}$  & {8}& {1.3e-1}& {3.1e-2} & {8}& {1.4e-1}& {3.5e-2} & {7}& {1.7e-1}& {4.0e-2} & {7}& {1.6e-1}& {4.3e-2} & {7}& {4.4e-2}& {2.4e-2} & {8}& {2.9e-2}& {1.6e-2} \\
				$2^{7}$  & {9}& {8.9e-2}& {1.6e-2} & {8}& {1.2e-1}& {2.2e-2} & {7}& {1.4e-1}& {2.6e-2} & {7}& {9.7e-2}& {2.3e-2} & {8}& {2.9e-2}& {1.2e-2} & {8}& {4.1e-2}& {1.2e-2} \\
				$2^{8}$  & {9}& {6.3e-2}& {8.5e-3} & {8}& {1.1e-1}& {1.3e-2} & {7}& {9.8e-2}& {1.7e-2} & {7}& {3.7e-2}& {1.2e-2} & {8}& {4.1e-2}& {1.0e-2} & {8}& {2.7e-2}& {6.3e-3} \\
				$2^{9}$  & {11}& {3.6e-2}& {4.1e-3} & {9}& {9.4e-2}& {8.3e-3} & {7}& {5.2e-2}& {1.2e-2} & {8}& {3.3e-2}& {9.7e-3} & {8}& {4.6e-2}& {1.1e-2} & {8}& {4.6e-2}& {7.9e-3} \\
				$2^{10}$  & {11}& {1.8e-2}& {2.2e-3} & {10}& {8.2e-2}& {5.1e-3} & {7}& {4.9e-2}& {1.1e-2} & {8}& {4.5e-2}& {1.1e-2} & {8}& {4.5e-2}& {9.6e-3} & {8}& {4.6e-2}& {7.7e-3} \\
				\hline
			\end{tabular}
		\end{center}
		\caption{Iterations to tolerance of $\mathcal{P}$-GMRES with related 2-norm relative error, and approximation error varying $N$ and the grid for  $\gamma=0.5$, $\beta=0.9$.}\label{tab_ex3}
	\end{table}

	\paragraph{Test 4} We have shown that $\mathcal{P}$-GMRES works in the case where $\gamma=0.5$ and that it outperforms PFCGS. Here we only focus on the behavior of $\mathcal{P}$-GMRES and show that it is a suitable {solver} even in the extreme anisotropic cases where $\gamma=0$ or $\gamma=1$.
	
	Table \ref{tab_ex4} shows It, the infinity norm error $e^\infty$ and the relative 2-norm error $e^\mathrm{rel}$ varying $\beta\in\lbrace 0.1,0.3,0.7\rbrace$ and $\gamma\in\lbrace 0,1\rbrace$. We recall that when \lq\lq -'' is displayed, the solver exceeded the maximum amount of iterations fixed to $100$. 
	
	When $\beta=0.1$ and $\beta=0.3$, despite the strong spatial anisotropy, the coefficient matrix is close to Hermitian (see Remark \ref{rem:laplacian}), and therefore $\mathcal{P}$-GMRES is expected to be a suitable preconditioner. Indeed, when considering $\mathbf{\epsilon}^{(1)},\mathbf{\epsilon}^{(4)},\mathbf{\epsilon}^{(6)}$, It does not increase with $N$ for both choices of $\gamma$ and both errors seem to decrease with order $1+\beta$ as observed in Test 2 with $\gamma=0.5$.
	
	When $\beta=0.7$, as the coefficient matrix is close to skew-symmetric, our multigrid preconditioning is not effective anymore. In fact, when $\gamma=0$ we observe stable iterations only for the composite meshes, but the error does not decrease as expected (see Test 2). When considering the grid mapped by $g_{q,\mathbf{\epsilon}}(x)$ with $\mathbf{\epsilon}=\mathbf{\epsilon}^{(6)}$, the iterations are not stable while increasing $N$, but the predicted convergence order is restored. When $\gamma=1$, $\mathcal{P}$-GMRES does not converge for any of the tested grids mapped by $g_{q,\mathbf{\epsilon}}(x)$. When using a composite mesh, instead, $\mathcal{P}$-GMRES converges but the number of iterations is large and increases with $N$, while the error decreases too slowly. 
	%
	%
	
	\begin{table}[!ht]\scriptsize\begin{center}
			\begin{tabular}{@{\extracolsep{3pt}}c@{\hspace{0.2cm}}|c@{\hspace{0.2cm}}c@{\hspace{0.2cm}}g@{\hspace{0.2cm}}c@{\hspace{0.1cm}}c@{\hspace{0.2cm}}g@{\hspace{0.2cm}}c@{\hspace{0.1cm}}c@{\hspace{0.2cm}}g@{\hspace{0.2cm}}c@{\hspace{0.1cm}}c@{\hspace{0.2cm}}g@{\hspace{0.2cm}}c@{\hspace{0.1cm}}c@{\hspace{0.2cm}}g@{\hspace{0.2cm}}c@{\hspace{0.1cm}}c@{\hspace{0.1cm}}c@{\hspace{0.1cm}}c@{\hspace{0.15cm}}c@{\hspace{0.1cm}}}
				\hline\vspace{-0.2cm} \\
				&&& \multicolumn{6}{c}{Composite mesh mapped by $g_G(N)$}&\multicolumn{9}{c}{Non-uniform mesh mapped by $g_{q,\mathbf{\epsilon}}(x)$}\\ \cline{4-9}\cline{10-18}\vspace{-0.2cm}\\
				\multirow{2}{*}{$\gamma$}&\multirow{2}{*}{$\beta$} &\multirow{2}{*}{$N\!+\!1$}  & \multicolumn{3}{c}{$g_G(N)=\lfloor \sqrt{N}\rfloor$} & \multicolumn{3}{c}{$g_G(N)=\lfloor \log_2{N}\rfloor$} & \multicolumn{3}{c}{$\mathbf{\epsilon}^{(1)}=(0.1,0.05)$} & \multicolumn{3}{c}{$\mathbf{\epsilon}^{(4)}=(0.45,0.05)$} & \multicolumn{3}{c}{$\mathbf{\epsilon}^{(6)}=(1,0)$}\\\cline{4-6}\cline{7-9}\cline{10-12}\cline{13-15}\cline{16-18}
				\vspace{-0.2cm} \\
				& & & It & $e^\infty$ & $e^{\text{rel}}$ & It & $e^\infty$ & $e^{\text{rel}}$ & It & $e^\infty$ & $e^{\text{rel}}$ & It & $e^\infty$ & $e^{\text{rel}}$ & It & $e^\infty$ & $e^{\text{rel}}$ \\ \hline\hline
				\vspace{-0.2cm} \\
				\multirow{20}{*}{0}&\multirow{6}{*}{0.1}
				& $2^{5}$  & {10}& {9.2e-4}& {1.3e-3} & {11}& {1.0e-3}& {1.4e-3} & {7}& {1.5e-3}& {2.0e-3} & {7}& {1.5e-3}& {2.0e-3} & {7}& {1.5e-3}& {2.1e-3} \\
				& & $2^{6}$  & {12}& {4.7e-4}& {6.2e-4} & {11}& {5.1e-4}& {6.7e-4} & {7}& {7.4e-4}& {9.5e-4} & {7}& {7.3e-4}& {9.9e-4} & {7}& {7.4e-4}& {1.0e-3} \\
				& & $2^{7}$  & {21}& {2.6e-4}& {3.3e-4} & {11}& {2.7e-4}& {3.4e-4} & {7}& {3.6e-4}& {4.6e-4} & {7}& {3.6e-4}& {4.8e-4} & {7}& {3.6e-4}& {4.9e-4} \\
				& & $2^{8}$  & {23}& {1.4e-4}& {1.8e-4} & {18}& {1.5e-4}& {1.8e-4} & {7}& {1.8e-4}& {2.2e-4} & {7}& {1.8e-4}& {2.3e-4} & {7}& {1.8e-4}& {2.4e-4} \\
				& & $2^{9}$  & {27}& {7.9e-5}& {9.5e-5} & {18}& {7.9e-5}& {9.5e-5} & {7}& {8.5e-5}& {1.1e-4} & {7}& {8.5e-5}& {1.1e-4} & {7}& {8.5e-5}& {1.1e-4} \\
				& & $2^{10}$  & {-}& {-}& {-} & {25}& {4.3e-5}& {5.1e-5} & {7}& {4.1e-5}& {5.1e-5} & {7}& {4.1e-5}& {5.4e-5} & {7}& {4.1e-5}& {5.5e-5} \\
				\cline{2-18}\vspace{-0.2cm} \\  
				&\multirow{6}{*}{0.3}  
				& $2^{5}$  & {14}& {6.1e-3}& {7.4e-3} & {14}& {7.4e-3}& {8.8e-3} & {17}& {6.5e-3}& {7.8e-3} & {17}& {4.5e-3}& {6.8e-3} & {17}& {4.6e-3}& {7.2e-3} \\
				& & $2^{6}$  & {12}& {3.1e-3}& {3.6e-3} & {14}& {3.9e-3}& {4.5e-3} & {13}& {2.7e-3}& {3.2e-3} & {13}& {2.0e-3}& {3.0e-3} & {13}& {2.0e-3}& {3.2e-3} \\
				& & $2^{7}$  & {12}& {1.8e-3}& {2.0e-3} & {11}& {2.2e-3}& {2.4e-3} & {9}& {1.1e-3}& {1.3e-3} & {10}& {9.0e-4}& {1.3e-3} & {10}& {9.0e-4}& {1.4e-3} \\
				& & $2^{8}$  & {20}& {1.1e-3}& {1.2e-3} & {12}& {1.3e-3}& {1.4e-3} & {10}& {4.5e-4}& {5.4e-4} & {11}& {3.9e-4}& {5.7e-4} & {11}& {3.9e-4}& {6.1e-4} \\
				& & $2^{9}$  & {16}& {6.9e-4}& {7.4e-4} & {15}& {7.5e-4}& {8.0e-4} & {11}& {1.9e-4}& {2.3e-4} & {10}& {1.7e-4}& {2.5e-4} & {10}& {1.7e-4}& {2.6e-4} \\
				& & $2^{10}$  & {30}& {4.3e-4}& {4.5e-4} & {20}& {4.5e-4}& {4.8e-4} & {10}& {7.9e-5}& {9.6e-5} & {10}& {7.2e-5}& {1.1e-4} & {10}& {7.2e-5}& {1.1e-4} \\
				\cline{2-18}\vspace{-0.2cm} \\  
				& \multirow{6}{*}{0.7}  
				& $2^{5}$  & {17}& {1.1e-1}& {1.2e-1} & {17}& {1.3e-1}& {1.4e-1} & {-}& {-}& {-} & {33}& {3.2e-2}& {3.7e-2} & {20}& {1.8e-2}& {3.2e-2} \\
				& & $2^{6}$  & {18}& {6.3e-2}& {6.9e-2} & {18}& {8.8e-2}& {9.4e-2} & {39}& {9.1e-2}& {9.5e-2} & {33}& {1.1e-2}& {1.3e-2} & {22}& {6.7e-3}& {1.2e-2} \\
				& & $2^{7}$  & {18}& {3.0e-2}& {3.1e-2} & {18}& {6.1e-2}& {6.3e-2} & {-}& {-}& {-} & {42}& {3.6e-3}& {4.3e-3} & {19}& {2.3e-3}& {4.3e-3} \\
				& & $2^{8}$  & {19}& {1.6e-2}& {1.7e-2} & {18}& {4.2e-2}& {4.3e-2} & {-}& {-}& {-} & {51}& {1.1e-3}& {1.4e-3} & {23}& {8.0e-4}& {1.5e-3} \\
				& & $2^{9}$  & {20}& {9.5e-3}& {9.6e-3} & {18}& {2.9e-2}& {3.0e-2} & {-}& {-}& {-} & {56}& {3.5e-4}& {4.6e-4} & {26}& {2.7e-4}& {5.1e-4} \\
				& & $2^{10}$  & {27}& {6.9e-3}& {7.0e-3} & {19}& {2.1e-2}& {2.1e-2} & {-}& {-}& {-} & {51}& {1.3e-4}& {1.9e-4} & {32}& {1.2e-4}& {2.3e-4} \\
				\hline\hline\vspace{-0.2cm}\\
				\multirow{20}{*}{1}&\multirow{6}{*}{0.1}
				& $2^{5}$  & {11}& {8.0e-5}& {6.1e-5} & {11}& {1.4e-4}& {8.6e-5} & {8}& {4.2e-4}& {3.3e-4} & {8}& {4.2e-4}& {3.6e-4} & {8}& {4.2e-4}& {3.7e-4} \\
				& & $2^{6}$  & {13}& {4.5e-5}& {5.4e-5} & {11}& {4.0e-5}& {3.8e-5} & {8}& {2.0e-4}& {1.3e-4} & {8}& {2.0e-4}& {1.4e-4} & {8}& {2.0e-4}& {1.4e-4} \\
				& & $2^{7}$  & {21}& {3.4e-5}& {3.9e-5} & {12}& {2.8e-5}& {3.2e-5} & {8}& {9.3e-5}& {4.8e-5} & {8}& {9.3e-5}& {5.3e-5} & {8}& {9.3e-5}& {5.4e-5} \\
				& & $2^{8}$  & {32}& {2.1e-5}& {2.4e-5} & {15}& {1.9e-5}& {2.2e-5} & {8}& {4.3e-5}& {1.7e-5} & {8}& {4.3e-5}& {1.9e-5} & {8}& {4.3e-5}& {1.9e-5} \\
				& & $2^{9}$  & {-}& {-}& {-} & {18}& {1.2e-5}& {1.4e-5} & {8}& {2.0e-5}& {5.9e-6} & {8}& {2.0e-5}& {6.3e-6} & {8}& {2.0e-5}& {6.3e-6} \\
				& & $2^{10}$  & {-}& {-}& {-} & {25}& {7.4e-6}& {8.0e-6} & {8}& {9.5e-6}& {2.2e-6} & {8}& {9.5e-6}& {2.1e-6} & {8}& {9.5e-6}& {2.1e-6} \\
				\cline{2-18}\vspace{-0.2cm} \\  
				& \multirow{6}{*}{0.3}  
				& $2^{5}$  & {16}& {1.3e-3}& {9.4e-4} & {19}& {1.1e-3}& {8.0e-4} & {21}& {1.3e-3}& {1.0e-3} & {19}& {4.6e-4}& {6.4e-4} & {21}& {4.5e-4}& {6.3e-4} \\
				& & $2^{6}$  & {15}& {8.9e-4}& {5.7e-4} & {15}& {8.0e-4}& {5.1e-4} & {13}& {5.6e-4}& {4.5e-4} & {15}& {2.2e-4}& {3.1e-4} & {15}& {2.2e-4}& {3.1e-4} \\
				& & $2^{7}$  & {20}& {5.7e-4}& {3.2e-4} & {13}& {5.2e-4}& {3.0e-4} & {12}& {2.3e-4}& {2.0e-4} & {10}& {1.0e-4}& {1.4e-4} & {12}& {1.0e-4}& {1.4e-4} \\
				& & $2^{8}$  & {26}& {3.5e-4}& {1.7e-4} & {15}& {3.3e-4}& {1.6e-4} & {11}& {9.6e-5}& {8.3e-5} & {11}& {4.7e-5}& {6.1e-5} & {11}& {4.7e-5}& {6.2e-5} \\
				& & $2^{9}$  & {-}& {-}& {-} & {16}& {2.1e-4}& {8.8e-5} & {11}& {3.9e-5}& {3.5e-5} & {11}& {2.1e-5}& {2.6e-5} & {11}& {2.1e-5}& {2.7e-5} \\
				& & $2^{10}$  & {-}& {-}& {-} & {21}& {1.3e-4}& {4.7e-5} & {11}& {1.6e-5}& {1.4e-5} & {11}& {8.8e-6}& {1.1e-5} & {11}& {8.8e-6}& {1.1e-5} \\
				\cline{2-18}\vspace{-0.2cm} \\  
				& \multirow{6}{*}{0.7}  
				& $2^{5}$  & {17}& {5.2e-2}& {1.4e-2} & {17}& {6.3e-2}& {1.6e-2} & {-}& {-}& {-} & {-}& {-}& {-} & {-}& {-}& {-} \\
				& & $2^{6}$  & {21}& {2.8e-2}& {5.2e-3} & {18}& {4.1e-2}& {7.5e-3} & {-}& {-}& {-} & {-}& {-}& {-} & {-}& {-}& {-} \\
				& & $2^{7}$  & {29}& {9.7e-3}& {1.5e-3} & {19}& {2.7e-2}& {3.5e-3} & {-}& {-}& {-} & {-}& {-}& {-} & {-}& {-}& {-} \\
				& & $2^{8}$  & {36}& {3.4e-3}& {6.0e-4} & {23}& {1.8e-2}& {1.6e-3} & {-}& {-}& {-} & {-}& {-}& {-} & {-}& {-}& {-} \\
				& & $2^{9}$  & {-}& {-}& {-} & {27}& {1.2e-2}& {7.7e-4} & {-}& {-}& {-} & {-}& {-}& {-} & {-}& {-}& {-} \\
				& & $2^{10}$  & {-}& {-}& {-} & {28}& {7.7e-3}& {3.7e-4} & {-}& {-}& {-} & {-}& {-}& {-} & {-}& {-}& {-} \\
				\hline
			\end{tabular}
		\end{center}
		\caption{Iterations of $\mathcal{P}$-GMRES with related 2-norm relative error, and approximation error varying $N$ and the grid for $\gamma\in\lbrace 0,1\rbrace$.}\label{tab_ex4}
	\end{table}
	
	\section{Conclusions}\label{sec:conclusion}
	We have considered a FVE discretization over a generic mesh of a conservative steady-state two-sided FDE whose solution is singular at the boundary, with a special focus on grids that combine a graded mesh near the singularity with a uniform mesh where the solution is smooth. The approximation order of the considered graded meshes has numerically been compared with the one of the composite mesh used in \cite{wang}, 
	showing that lower approximation errors can be obtained.
	
	
	By exploiting the mesh structure, we have computed the symbol of the resulting coefficient matrices in case of non-uniform meshes mapped by a function. The related spectral information has been leveraged to build an ad-hoc multigrid preconditioner. Through a wide number of numerical tests, we have shown that, except for $\gamma\approx0$ or $\gamma\approx1$ and $\beta\approx1$, such a multigrid is a valid alternative to the circulant preconditioner developed in \cite{wang} for composite meshes.

	
	
	
	When $\beta\approx 0$, as done in \cite{paper_MGS}, a band approximation of the coefficient matrix could be exploited in order to further reduce the computational cost of the preconditioning iteration. We stress that multigrid methods should perform even better in the two-dimensional case with respect to the circulant preconditioner since it is well-known that {multilevel circulant matrices} used as preconditioners for multilevel Toeplitz matrices cannot ensure superlinear convergence; see~\cite{circulant}.
	
	Finally, we mention that all the retrieved spectral results could easily be extended to time-dependent problems treated, e.g., in \cite{simmons} and this could be used to have insights on the stability of the chosen time scheme.
	
}

\section*{Acknowledgments} 
This research was supported by the GNCS-INDAM (Italy), the Swiss National Science Foundation SNF via the projects Stress-Based Methods for Variational Inequalities in Solid Mechanics n. 186407 and ExaSolvers n. 162199, and the EuroHPC TIME-X project.


%
%




\appendix

\section{Proof of Proposition \ref{prop_existence_grid}}\label{appendix_2}
The explicit form of coefficients $a,b,c,m,q$ is obtained by solving the following equation
\begin{equation}\label{eq_temp_ls}
\begin{cases}
g_{q,\mathbf{\epsilon}}(\epsilon_1^-)&=g_{q,\mathbf{\epsilon}}(\epsilon_1^+)\\
g_{q,\mathbf{\epsilon}}((\epsilon_1+\epsilon_2)^-)&=g_{q,\mathbf{\epsilon}}((\epsilon_1+\epsilon_2)^+)\\
g_{q,\mathbf{\epsilon}}'(\epsilon_1^-)&=g_{q,\mathbf{\epsilon}}'(\epsilon_1^+)\\
g_{q,\mathbf{\epsilon}}'((\epsilon_1+\epsilon_2)^-)&=g_{q,\mathbf{\epsilon}}'((\epsilon_1+\epsilon_2)^+)\\
g_{q,\mathbf{\epsilon}}(1)&=1
\end{cases}
\end{equation}
where $g_{q,\mathbf{\epsilon}}(\xi^\pm)=\lim_{x\to\xi^\pm}g_{q,\mathbf{\epsilon}}(x)$. Equation \eqref{eq_temp_ls} can be seen as a linear system with coefficient matrix 
\begin{equation*}
G=\begin{pmatrix}
\epsilon_1^2    &    \epsilon_1 &  1   &  0  &     0 \\
(\epsilon_1+\epsilon_2)^2  &  \epsilon_1+\epsilon_2 &  1 & -(\epsilon_1+\epsilon_2) & -1\\
2\epsilon_1  &    1  &  0  &   0    &   0 \\
2(\epsilon_1+\epsilon_2) &  1&    0   &  -1  &    0 \\
0      &  0  &  0 &    1   &    1 
\end{pmatrix},
\end{equation*}
whose determinant is $\mathrm{det}(G)=2\epsilon_1\epsilon_2-2\epsilon_2+\epsilon_2^2$. Finally, since $\mathrm{det}(G)=0$ if and only if $\epsilon_1=\frac{2-\epsilon_2}{2}$, and $0<\epsilon_1+\epsilon_2\leq 1$ with $\epsilon_2>0$, we conclude that $\mathrm{det}(G)\neq0$ and  therefore $g_{q,\mathbf{\epsilon}}(x)$ is well-defined.

\section{Proof of Theorem \ref{th_symbol_val_sing}}\label{appendix_1}
Let us fix $N\in\mathbb{N}$ and let $\lbrace \hat{x}_i\rbrace_{i=0}^{N+1}$ be the uniform grid with $\hat{x}_i=ih$, $i=0,...,N+1$, and $h=\tfrac{1}{N+1}$. Then, {according to equation \eqref{eq_graded_mesh}, by} letting $x_i=g(\hat{x}_i),\ i=0,...,N+1$, from the Taylor expansion of $g(\hat{x}_{i-1})$, $\forall i$, it holds
\begin{equation*}
g(\hat{x}_{i-1})=g(\hat{x}_{i})-g'(\hat{x}_{i})h+g''(\hat{x}_{i})\frac{h^2}{2}+\O(h^3),\qquad \forall i=1,...,N,
\end{equation*}
and since $h_i=g(\hat{x}_{i})-g(\hat{x}_{i-1})$, we have
\begin{equation}\label{eq_hi}
h_i=g'(\hat{x}_{i})h-g''(\hat{x}_{i})\frac{h^2}{2}+\O(h^3).
\end{equation}
In the case where $k\in\mathbb{Z}$, from \eqref{eq_hi} we have
\begin{equation}\label{eq_hik}
h_{i+k}=g'(\hat{x}_{i+ k})h-g''(\hat{x}_{i+ k})\frac{h^2}{2}+\O(h^3),
\end{equation}
and through the Taylor expansions of $g'(\hat{x}_{i+k})$ and $g''(\hat{x}_{i+k})$ we finally obtain
\begin{align*}
h_{i+k}=g'(\hat{x}_{i})h+g''(\hat{x}_{i})h^2\frac{2k-1}{2}+\O(k^2h^3).
\end{align*}
Then, with $\gamma=\frac{1}{2}$, $K(x)=K$ and $\tilde{K}=\frac{K}{\Gamma(\beta+1)}$, from equation \eqref{eq_LINEAR_system_discretizedeq} we have
\begin{equation*}
a_{i,i}=\tilde{K}\left[
\frac{(\frac{h_i}{2})^\beta}{h_i}+\frac{1}{2}\frac{(\frac{h_{i}}{2})^\beta-(h_{i+1}+\frac{h_{i}}{2})^\beta}{h_{i+1}}
-\frac{1}{2}\frac{(h_i+\frac{h_{i+1}}{2})^\beta-(\frac{h_{i+1}}{2})^\beta}{h_i}+\frac{(\frac{h_{i+1}}{2})^\beta}{h_{i+1}}
\right]
=\tilde{K}\left(S_1+\frac{1}{2}S_2-\frac{1}{2}S_3+S_4\right),
\end{equation*}
where
\begin{align*}
S_1&=\frac{(\frac{h_i}{2})^\beta}{h_i}
=\frac{1}{2^\beta}h_i^{\beta-1}
=\frac{1}{2^\beta}\left(\gi h+\O(h^2)\right)^{\beta-1}
=\frac{1}{2^\beta(\gi h)^{1-\beta}}\left(1+\O(h)\right);\\
%
S_2&=\frac{(\frac{h_i}{2})^\beta-(h_{i+1}+\frac{h_i}{2})^\beta}{h_{i+1}}
=\frac{(\gi h)^\beta\left(\big(\frac{1}{2}+\O(h)\big)^\beta-\big(\frac{3}{2}+\O(h)\big)^\beta\right)}{\gi h(1+\O(h))}
=\frac{1}{2^\beta(\gi h)^{1-\beta}}\frac{\big(1+\O(h)\big)^\beta-3^\beta\big(1+\O(h)\big)^\beta}{1+\O(h)}\\
&
=\frac{1}{2^\beta(\gi h)^{1-\beta}} \left(1-3^\beta+\O(h)\right) \Big(1+\O(h)\Big)
=\frac{1}{2^\beta(\gi h)^{1-\beta}}\left(1-3^\beta+\O(h)\right);\\
S_3&=\frac{(h_{i}+\frac{h_{i+1}}{2})^\beta-(\frac{h_{i+1}}{2})^\beta}{h_{i}}
=\frac{(\gi h)^\beta\left(\big(\frac{3}{2}+\O(h)\big)^\beta-\big(\frac{1}{2}+\O(h)\big)^\beta\right)}{\gi h(1+\O(h))}
=\frac{1}{2^\beta(\gi h)^{1-\beta}} \frac{3^\beta\big(1+\O(h)\big)^\beta-\big(1+\O(h)\big)^\beta}{1+\O(h)}\\
&
=\frac{1}{2^\beta(\gi h)^{1-\beta}} \left(3^\beta-1+\O(h)\right) \Big(1+\O(h)\Big)
=\frac{1}{2^\beta(\gi h)^{1-\beta}}\left(3^\beta-1+\O(h)\right);\\
S_4&=\frac{(\frac{h_{i+1}}{2})^\beta}{h_{i+1}}
=\frac{1}{2^\beta}h_{i+1}^{\beta-1}
=\frac{1}{2^\beta}\left(\gi h+\O(h^2)\right)^{\beta-1}
=\frac{1}{2^\beta(\gi h)^{1-\beta}}\left(1+\O(h)\right);
\end{align*}
with $g'_i=\g$ and $g''_i=g''(\hat{x}_i)$.\\
Assembling $a_{i,i}$ we obtain
\begin{align}\label{eq_aii_approx}
h^{1-\beta}a_{i,i}=\frac{\tilde{K}}{2^\beta{\gi}^{1-\beta}}\left(3-3^\beta\right)+\mathrm{O}(h).
\end{align}
With the same approach we obtain
\begin{align}                                
h^{1-\beta}a_{i,i\pm 1}=&\frac{\tilde{K}}{2^\beta {\gi}^{1-\beta}}\left[
3^{\beta+1}-4-5^\beta
\right]+\mathrm{O}(h);\label{eq_ai1_approx}\\
h^{1-\beta}a_{i,i\pm k}=&\frac{\tilde{K}}{2^\beta{\gi}^{1-\beta}}\left[
3(2k+1)^\beta-3(2k-1)^\beta+(2k-3)^\beta-(2k+3)^\beta
\right]+\mathrm{O}(kh),\label{eq_aik_approx}
\end{align}
with $1<k\leq N^{\q},\ 0<\q<1,$ such that $kh\rightarrow 0$ as $k\rightarrow\infty$.\\ 

In the case where $N^{\q}<k\leq N$ the approximation yields a large error, therefore we prove that $a_{i,i+ k}=\mathrm{o}(h)$. Let $r=\sum_{j=1}^{k}h_{i+j}$, then $0<r<1\ \forall k$ and
\begin{itemize}
	\item if $k=\O(N^\q)=\O(h^{-\q})$, we have
	\begin{equation}\label{eq_r}
	r=\sum\limits_{j=1}^{\O(N^\q)}h_{i+j}
	=\sum\limits_{j=1}^{\O(N^\q)}\Big(\gi h+\O(jh^2)\Big)
	=\gi \O(N^\q)h+\O(N^{2\q}h^2)=\gi\O(h^{1-\q}),
	\end{equation}
	\item while if $k=\O(N)=\O(h^{-1})$, $r$ is a constant independent of $N$.
\end{itemize}
\noindent By collecting $r$, in $a_{i,i+ k}$ we have terms of the form $(1+\tilde{h})^\beta$, with $\tilde{h}=\frac{h_i}{2r},\frac{h_{i}}{2r}-\frac{h_{i+k}}{r},\ldots \ $. In order to use the Taylor expansion of $(1+\tilde{h})^\beta$ we first need to prove that $\tilde{h}\rightarrow 0$ as $N\rightarrow\infty$. We divide the analysis in two cases:
\begin{itemize}
	\item[1)] if $g'(\hat{x})\neq 0$ in $[0,1]$, then from equations \eqref{eq_hi} and \eqref{eq_r} we have
	\begin{align*}
	\frac{h_{i+k}}{r}=\frac{g_{i+k}' h+\O(h^2)}{\gi\O(h^{1-\q})}=\O(h^\q),
	\end{align*}
	which tends to zero as $N\rightarrow\infty$ for any $0<\q<1$.
	\item[2)] if $g'(\hat{x})$ vanishes in $[0,1]$, the worst possible scenario happens when $h_i\gg h_{i+k}$. Without restrictions to the general case we assume that $g'(\hat{x})$ has a zero of order $\ord$ at $\hat{x}_0=0$, hence $g(\hat{x})\approx\hat{x}^{\ord+1}$ when $\hat{x}\rightarrow 0$. Then by considering $k=-N^\q$ and $i=N^\q+1$, such that $h_{N^\q+1}=h_i\gg h_{i+k}=h_1$, we have
	{\small\begin{align*}
	\frac{h_{N^\q+1}}{\sum_{j=1}^{N^\q}h_{N^\q+1-j}}
	&=\frac{x_{N^\q+1}-x_{N^\q}}{h_1+h_2+...+h_{N^\q}}
	=\frac{g(\hat{x}_{N^\q+1})-g(\hat{x}_{N^\q})}{g(\hat{x}_{N^\q})-g(\hat{x}_{0})}
	=\frac{\left(\frac{N^\q+1}{N}\right)^{\ord+1}-\left(\frac{N^\q}{N}\right)^{\ord+1}}{\left(\frac{N^\q}{N}\right)^{\ord+1}}
	=\left(1+\frac{1}{N^\q}\right)^{\ord+1}-1=\O(h^{\q\ord}),
	\end{align*}}
	which tends to zero as $N\rightarrow\infty$ for any $0<\q<1$.
\end{itemize}
\noindent Now, let $N^\q<k\leq N$ and approximate the coefficient $a_{i,i+ k}$ as follows:
\begin{align*}
a_{i,i+ k}=
&\frac{\tilde{K}}{2}\Bigg[\frac{r^\beta(1+\frac{h_i}{2r})^\beta-r^\beta(1+\frac{h_i}{2r}-\frac{h_{i+k}}{r})^\beta}{h_{i+k}}+
\frac{r^\beta(1+\frac{h_i}{2r})^\beta-r^\beta(1+\frac{h_i}{2r}+\frac{h_{i+k+1}}{r})^\beta}{h_{i+k+1}}+\\
&-\frac{r^\beta(1-\frac{h_{i+1}}{2r})^\beta-r^\beta(1-\frac{h_{i+1}}{2r}-\frac{h_{i+k}}{r})^\beta}{h_{i+k}}-
\frac{r^\beta(1-\frac{h_{i+1}}{2r})^\beta-r^\beta(1-\frac{h_{i+1}}{2r}+\frac{h_{i+k}}{r})^\beta}{h_{i+k+1}}\Bigg]\\
=&\frac{\tilde{K} r^\beta}{2h_{i+k}h_{i+k+1}}\Bigg[h_{i+k+1}\left((1+\frac{h_i}{2r})^\beta-(1+\frac{h_i}{2r}-\frac{h_{i+k}}{r})^\beta-(1-\frac{h_{i+1}}{2r})^\beta+(1-\frac{h_{i+1}}{2r}-\frac{h_{i+k}}{r})^\beta\right)+\\
&+h_{i+k}\left((1+\frac{h_i}{2r})^\beta-(1+\frac{h_i}{2r}+\frac{h_{i+k+1}}{r})^\beta-(1-\frac{h_{i+1}}{2r})^\beta+(1-\frac{h_{i+1}}{2r}+\frac{h_{i+k}}{r})^\beta\right)\Bigg].
\end{align*}
By replacing each $(1+ \tilde{h} )^\beta$ with its Taylor expansion 
$$(1+\tilde{h})^\beta=1+\beta \tilde{h}+ \frac{\beta(\beta-1)}{2}\tilde{h}^2+\mathrm{O}(\tilde{h}^3),$$
we observe an exact cancellation of the terms of degree $0$ and $1$ inside the square brackets in $a_{i,i+ k}$. The exact cancellation happens even for the term of degree $2$ but it is harder to see, therefore we report the computations below:
\begin{align}
a_{i,i+ k}=&\frac{\tilde{K} r^\beta \beta(\beta-1)}{4 h_{i+k}h_{i+k+1}} 
\Bigg[
h_{i+k+1}\bigg(\frac{h_i^2}{4r^2}+\mathrm{O}(h_i^3)-(\frac{h_i}{2r}-\frac{h_{i+k}}{r})^2+\mathrm{O}((h_i-h_{i+k})^3)-\frac{h_{i+1}^2}{4r^2}+\mathrm{O}(h_{i+1}^3)+\nonumber\\
&\qquad\qquad\qquad\qquad\quad+(\frac{h_{i+1}}{2r}+\frac{h_{i+k}}{r})^2+\mathrm{O}((h_{i+1}+h_{i+k})^3))\bigg)+\nonumber\\
&\qquad\qquad\quad\ \  +h_{i+k}\bigg(\frac{h_i^2}{4r^2}+\mathrm{O}(h_i^3)-(\frac{h_i}{2r}+\frac{h_{i+k+1}}{r})^2+\mathrm{O}((h_i+h_{i+k+1})^3)-\frac{h_{i+1}^2}{4r^2}+\mathrm{O}(h_{i+1}^3)+\nonumber\\
&\qquad\qquad\qquad\qquad\quad+(-\frac{h_{i+1}}{2r}+\frac{h_{i+k+1}}{r})^2+\mathrm{O}((-h_{i+1}+h_{i+k+1})^3))\bigg)
\Bigg]\nonumber\\
=&\frac{\tilde{K} r^{\beta-2} \beta(\beta-1)}{4 h_{i+k}h_{i+k+1}} \left( h_{i+k+1}(h_i h_{i+k}+h_{i+1}h_{i+k})-h_{i+k}(h_i h_{i+k+1}+h_{i+1}h_{i+k+1})+\mathrm{O}(h^4)\right)\nonumber\\
=&\frac{r^{\beta-2}}{{g'_{i+k}}^2}\mathrm{O}(h^2).\nonumber
\end{align}
In the case where 
\begin{itemize}
	\item $k=\O(N^\q)$, from equation \eqref{eq_r} we have
	\begin{equation}\label{eq_aik_approx_k_large}
	h^{1-\beta}a_{i,i+ k}=
	\mathrm{O}(h^{1-\beta+2+(1-\q)(\beta-2)})=\mathrm{O}(h^{1-\beta+2+\beta-2+2\q-\beta \q})=\mathrm{o}(h^{1+\q}),
	\end{equation}
	\item $k=\O(N)$, since $r$ has a constant value, we have
	\begin{equation}\label{eq_aik_approx_k_large2}
	h^{1-\beta}a_{i,i+ k}=\mathrm{O}(h^{1-\beta+2})=\mathrm{o}(h^{2}).
	\end{equation}
\end{itemize}
Let $B_{N,M}$ be a diagonal-times-Toeplitz banded matrix of the form 
$$B_{N,M}=D_N(d(x))T_N(p^{\beta}_M(\theta)),$$
with $d(x)= \frac{\tilde{K}}{2^\beta (g'(x))^{1-\beta}}$ and $p^{\beta}_M(\theta)$ being the symbol in equation \eqref{eq_symbol}. {From Proposition \ref{prop_distribution} we have}
\begin{equation}\label{eq_symbol_B}
\lbrace\lbrace B_{N,M}\rbrace_N\rbrace_M\sim_\sigma \left( d(x)p^{\beta}_M(\theta), [0,1]\times[-\pi,\pi]\right).
\end{equation}
We now prove that $\lbrace\lbrace B_{N,M}\rbrace_N\rbrace_M$ is an a.c.s for $\lbrace h^{1-\beta}A_{N}\rbrace_N$.

Suppose that $g'(\tilde{x})\neq0$ in $[0,1]$, then, by choosing $M=N^\q$, from equations \eqref{eq_aii_approx}, \eqref{eq_ai1_approx}, \eqref{eq_aik_approx} and \eqref{eq_aik_approx_k_large} we have that matrix $h^{1-\beta}A_N-B_{N,N^\q}$ is a ``symmetric Toeplitz'' matrix, whose first row is 
\begin{equation}\label{eq_approx_th_sing_val}
\begin{matrix}
\O(h)\quad \ & \O(h)\quad \  & \O(2h)\quad \  
& \O(3h)\quad \  & \cdots \quad \  & \O(N^{\q}h)\quad \  
& \mathrm{o}(h^{1+\q})\quad \  & \cdots \quad \  &  \mathrm{o}(h^{1+\q})  \vspace{-10pt}\\
\ & \multicolumn{5}{l}{\begin{matrix}
	\underbrace{\qquad\ \qquad\ \qquad\ \qquad\ \qquad\ \qquad\ \qquad\ \qquad\ \quad\ }\\ N^{\q}\ \text{coefficients}
	\end{matrix}}
& \multicolumn{3}{c}{\begin{matrix}
	\underbrace{\qquad\ \qquad\ \qquad\ \qquad\ \quad\ }\\ N-N^\q-1\ \text{coefficients}
	\end{matrix}} \quad . 
\end{matrix}
\end{equation}
Note that the ``symmetric Toeplitz'' structure holds only while keeping $\O(\cdot)$ and $\mathrm{o}(\cdot)$. When we replace $\O(\cdot)$ and $\mathrm{o}(\cdot)$ with the exact values the structure could be lost.\\
Thanks to the structure we have 
\begin{equation*}
\begin{split}
\norm{h^{1-\beta}A_N-B_{N,N^\q}}_1 &\leq2\left(\O(h)+\sum_{k=1}^{N^\q}\O(kh)+\mathrm{o}(h^{1+\q})(N-N^\q-1)\right)\\
&=\O(h)+\O(N^{2\q}h)+\mathrm{o}((N-N^\q-1)h^{1+\q})\\
&=\O(h)+\O(h^{1-2\q})+\mathrm{o}(h^{\q}-h-h^{1+\q}),\\
\norm{h^{1-\beta}A_N-B_{N,N^\q}}_\infty &\leq \O(h)+\O(h^{1-2\q})+\mathrm{o}(h^{\q}-h-h^{1+\q}),
\end{split}
\end{equation*}
and through the Hölder inequality,
\begin{equation}\label{eq_Norm_N}
\begin{split}
\norm{h^{1-\beta}A_N-B_{N,N^\q}}_2\leq \sqrt{\norm{h^{1-\beta}A_N-B_{N,N^\q}}_1\norm{h^{1-\beta}A_N-B_{N,N^\q}}_\infty}
&\leq \O(h)+\O(h^{1-2\q})+\mathrm{o}(h^{\q}-h-h^{1+\q}),
\end{split}
\end{equation}
which tends to zero as $N\rightarrow\infty$ if $0<\q<\frac{1}{2}$. From Definition \ref{def_acs} it follows that $\lbrace\lbrace B_{N,M}\rbrace_{N}\rbrace_M$ is an a.c.s for \cred{$\lbrace h^{1-\beta}A_N\rbrace_N$}, and from equation \eqref{eq_symbol_B} and Theorem \ref{th_simbolo_acs} we have the thesis, {since from Proposition \ref{prop_spectral_info} it holds that $p^{\beta}_M(\theta)$ converges to $p^{\beta}(\theta)$.}

Suppose now that there exist $\tilde{x}^{(1)},...,\tilde{x}^{(s)}\in[0,1]$ such that $g'(\tilde{x}^{(k)})=0,\ \forall k$ and consider the intervals $B(\tilde{x}^{(k)},\frac{1}{M})=\lbrace \tilde{x}\in[0,1]\ :\ \abs{\tilde{x}-\tilde{x}^{(k)}}<\frac{1}{M}\rbrace$. The function $g'(x)$ is continuous and strictly positive on $[0,1]\setminus I_M\ \forall M$, where
$$I_M=\bigcup_{k=1}^s B\Big(\tilde{x}^{(k)},\frac{1}{M}\Big),$$
therefore we write $h^{1-\beta}A_N-B_{N,M}=N_{N,M}+R_{N,M}$, where matrices $N_{N,M},R_{N,M}$ have small-norm and low-rank, respectively.\\
If we define $\tilde{a}_{i,j}=\left(h^{1-\beta}A_N-B_{N,M}\right)_{i,j}$, the matrix $R_{N,M}$ is sparse and its entries are
\begin{equation*}
r_{i,j}=
\begin{cases}
\tilde{a}_{i,j},& \text{if}\ \hat{x}_i\in I_M\ \text{or}\ \hat{x}_j\in I_M\\
0, & \text{otherwise}.
\end{cases}
\end{equation*}
Hence, 
given the equispaced grid $\hat{x}_i=ih,\ i=1,...,N$, then
\begin{equation*}
\text{rank}\left( R_{N,M} \right)\leq 2\vert I_M\vert
\leq 2s\Big(\frac{2/M}{h}+1\Big)
=2s\left(\frac{2}{M}+\frac{2}{NM}+\frac{1}{N}\right)N.
\end{equation*}
Moreover, as $R_{N,M}$ contains every coefficient $\tilde{a}_{i,j}$ which has a $g'_i$, with $i\in I_M$ in the denominator, the matrix $N_{N,M}$ can be approximated as in \eqref{eq_Norm_N}. In conclusion, for each $M$ there exists $N_M$ such that for $N>N_M$, rank$(R_{N,M})\leq \frac{3sN}{M}$ and $\norm{N_{N,M}}_\infty\leq \frac{1}{M}$, which, from Definition \ref{def_acs}, means that $\lbrace\lbrace B_{N,M}\rbrace_{N}\rbrace_M$ is an a.c.s for $\lbrace h^{1-\beta}A_{N}\rbrace_{N}$ and the thesis again follows from equation \eqref{eq_symbol_B} and Theorem \ref{th_simbolo_acs}.

\section{Proof of Theorem \ref{th_symbol}}\label{appendix_3}

The thesis is proven by combining Theorem \ref{th_symbol_val_sing} with \textbf{GLT1} and \textbf{GLT5}. Therefore, we only need to show that \textbf{GLT5} holds for the matrix sequence in equation \eqref{eq_temp_th_valsing}. We recall that \textbf{GLT5} consists in proving that
$$\mathrm{lim}_{N\to\infty}\frac{\norm{h^{1-\beta}A_N-h^{1-\beta}A_N^{\mathrm{H}}}_{\mathrm{tr}}}{N}{=0},$$
with $h=\frac{1}{N+1}$.\\
Let us now denote $\tilde{a}_{i,j}=h^{1-\beta}\left(A_N-A_N^{\mathrm{H}}\right)_{ij}$. Then, since equations \eqref{eq_aik_approx_k_large} and \eqref{eq_aik_approx_k_large2} hold also for negative values of $k$, we have 
$$\tilde{a}_{i,i+k}=h^{1-\beta}({a}_{i,i+k}-{a}_{i+k,i})=h^{1-\beta}({a}_{i,i+k}-{a}_{j,j-k})=O(h^{1+\q}),$$
with $j=i+k$, for $N^\q< \abs{k}\leq N$. From equation \eqref{eq_aik_approx} we have,
\begin{align*}
\tilde{a}_{i,i+k}=h^{1-\beta}({a}_{i,i+k}-{a}_{j,j-k})=&
\frac{\tilde{K}}{2^\beta{\gi}^{1-\beta}}\left[3(2k+1)^\beta-3(2k-1)^\beta+(2k-3)^\beta-(2k+3)^\beta\right]+\mathrm{O}(kh)+\\
&-\left(\frac{\tilde{K}}{2^\beta{g_j'}^{1-\beta}}\left[3(2k+1)^\beta-3(2k-1)^\beta+(2k-3)^\beta-(2k+3)^\beta\right]+\mathrm{O}(kh)\right)\\
=&\left(\frac{1}{{\gi}^{1-\beta}}-\frac{1}{g_{i+k}'^{1-\beta}}\right)\frac{\tilde{K}}{2^\beta}\left[3(2k+1)^\beta-3(2k-1)^\beta+(2k-3)^\beta-(2k+3)^\beta\right]+\mathrm{O}(kh)\\
=&\ O(kh),
\end{align*}
for $0\leq \abs{k} \leq N^\q$.
Therefore, a similar reasoning to the one made for $h^{1-\beta}A_N-B_{N,N^\q}$ in equation \eqref{eq_approx_th_sing_val} can be done also for $h^{1-\beta}\left(A_N-A_N^{\mathrm{H}}\right)$. Finally, from H\"older inequality it follows
\begin{equation*}
\frac{h^{1-\beta}\norm{A_N-A_N^{\mathrm{H}}}_{\mathrm{tr}}}{N}\leq \frac{Nh^{1-\beta}\norm{A_N-A_N^{\mathrm{H}}}_{2}}{N}\leq \O(h)+\O(h^{1-2\q})+\mathrm{o}(h^{\q}-h-h^{1+\q}),
\end{equation*}
which tends to $0$ as $N\to\infty$ and this concludes the proof.

\end{document}